\documentclass[smallextended]{article}
\usepackage[english]{babel}
\usepackage{graphicx}
\usepackage{amsthm}
\usepackage{exscale}
\usepackage{amsfonts}
\usepackage{amssymb}
\usepackage[intlimits,sumlimits]{amsmath}
\usepackage{amsmath,cases}
\usepackage[latin1]{inputenc}
\usepackage{mathrsfs}
\usepackage{fancyhdr}
\usepackage[misc,geometry]{ifsym}
\usepackage{color}
\usepackage{bm}
\usepackage{tikz}
\usepackage{multicol}
\let\svtikzpicture\tikzpicture
\def\tikzpicture{\noindent\svtikzpicture}
\usetikzlibrary{shapes,arrows}

\definecolor {bianco}{RGB}{255,255,255}

\tikzstyle{block} = [rectangle, fill=bianco,
    text width=11em, text centered, rounded corners, minimum height=4em]
\tikzstyle{line} = [draw, -latex']

\theoremstyle{definition}

\numberwithin{equation}{section}

\newtheorem{thm}{Theorem}[section]

\newtheorem{prop}[thm]{Proposition}
\newtheorem{lm}[thm]{Lemma}

\newtheorem{remark}[thm]{Remark}

\usepackage{a4wide}
\newcommand{\ud}{\mathrm{d}}
\newcommand{\mb}{\mathbf{m}}
\newcommand{\xb}{\mathbf{x}}
\newcommand{\yb}{\mathbf{y}}

\newcommand{\tb}{\mathbf{t}}

\newcommand{\thetab}{\boldsymbol{\theta}}
\newcommand{\etab}{\boldsymbol{\eta}}

\newcommand{\xisun}{\tilde{\boldsymbol{\xi}}^{(n)}}

\newcommand{\xinsu}{\tilde{\boldsymbol{\xi}}^{(n)}}
\newcommand{\ind}{1 \! \textrm{l}}

\def\rd{\mathbb{R}^d}
\def\rk{\mathbb{R}^k}
\def\naturals{\mathbb{N}}
\def\reals{\mathbb{R}}

\newcommand{\pp}{\textsf{P}}
\newcommand{\ee}{\textsf{E}}
\newcommand{\probabilityspace}{(\Omega, \mathscr{F}, \textsf{P})}

\newcommand{\pms}{[\Space]}

\newcommand{\randommeasure}{\tilde{\mathfrak{p}}}

\newcommand{\empiric}{\tilde{\mathfrak{e}}_n}
\newcommand{\empiricm}{\tilde{\mathfrak{e}}_{n,m}}
\newcommand{\pfrak}{\mathfrak{p}}
\newcommand{\xitil}{\tilde{\xi}}
\newcommand{\model}{\mathcal{M}}

\newcommand{\Space}{\mathbb{S}}

\author{Emanuele Dolera $\!\!^{\ast}$ and Eugenio Regazzini \\ \emph{Universit\`a degli Studi di Pavia}}

\title{\textbf{Uniform rates of the Glivenko-Cantelli convergence and their use in approximating Bayesian inferences}}
\date{}

\begin{document}
\maketitle

\begin{abstract}
This paper deals with suitable quantifications in approximating a probability measure by an ``empirical'' random probability measure $\hat{\pfrak}_n$, depending on the first $n$ terms of a sequence $\{\xitil_i\}_{i \geq 1}$ of random elements. Section \ref{sect:iid} studies the range of oscillation near zero of the Wasserstein distance $\ud^{(p)}_{\pms}$ between $\pfrak_0$ and $\hat{\pfrak}_n$, assuming the $\xitil_i$'s i.i.d.\! from $\pfrak_0$. In Theorem \ref{thm:fournier} $\pfrak_0$ can be fixed in the space of all probability measures on $(\rd, \mathscr{B}(\rd))$ and $\hat{\pfrak}_n$ coincides with the empirical measure $\empiric := \frac{1}{n} \sum_{i=1}^n \delta_{\xitil_i}$. In Theorem \ref{thm:Gauss} (Theorem \ref{thm:ExpFam}, respectively) $\pfrak_0$ is a $d$-dimensional Gaussian distribution (an element of a distinguished statistical exponential family, respectively) and $\hat{\pfrak}_n$ is another $d$-dimensional Gaussian distribution with estimated mean and covariance matrix (another element of the same family with an estimated parameter, respectively). These new results improve on allied recent works by providing 
also uniform bounds with respect to $n$, meaning the finiteness of the $p$-moment of $\sup_{n \geq 1} b_n \ud^{(p)}_{\pms}(\pfrak_0, \hat{\pfrak}_n)$ is proved for some diverging sequence $b_n$ of positive numbers. In Section  \ref{sect:exchangeable}, assuming the $\xitil_i$'s exchangeable, one studies the range of oscillation near zero of the Wasserstein distance between the conditional distribution---also called posterior---of the directing measure of the sequence, given $\xitil_1, \dots, \xitil_n$, and the point mass at $\hat{\pfrak}_n$. Similarly, a bound for the approximation of predictive distributions is given. Finally, Theorems from \ref{thm:BayesNP} to \ref{thm:BayesExpFam} reconsider Theorems from \ref{thm:fournier} to \ref{thm:ExpFam}, respectively, according to a Bayesian perspective.



\end{abstract}

\section{Introduction and description of the results} \label{sect:intro}

A recurrent problem in different branches of science concerns the approximation of a probability measure (p.m., for short)---either fixed or random---by means of another random p.m., say $\hat{\pfrak}_n$, of an empirical nature, as function of a random sample $\xinsu := (\xitil_1, \dots, \xitil_n)$ from the original measure. The accuracy of such approximation is usually measured in terms of a distinguished form of distance between p.m.'s, and the final objective is to study the stochastic convergence to zero of that distance between the original p.m. and $\hat{\pfrak}_n$, as $n$ goes to infinity.

This problem arose in connection with the issue, of a statistical nature, of estimating an unknown probability distribution (p.d., for short) on $\reals$ by means of the empirical frequency distribution. Thus, the original formulation considered $\xinsu$ as the initial $n$-segment of a sequence $\{\xitil_i\}_{i \geq 1}$ of i.i.d.\! real random variables with common p.d.\!  $\pfrak_0$, and $\hat{\pfrak}_n$ equal to $\empiric := \frac{1}{n} \sum_{i=1}^n \delta_{\xitil_i}$. In the same year and journal, Cantelli \cite{cantelli}, Glivenko \cite{glivenko} and Kolmogorov \cite{kolmoGC} published their respective solutions concerning the almost sure convergence to zero of the so-called uniform (or Kolmogorov) distance between $\pfrak_0$ and $\empiric$, while de Finetti \cite{definetti33} gave a solution in terms of the L\'{e}vy distance. Since this achievement was considered by some authors 
so important to earn the title of ``fundamental theorem of mathematical statistics'', a flourishing line of research developed from it with a view to providing refinements and extensions in various directions, such as: i) replacing $\empiric$ with smoothed versions of it; ii) providing central limit theorems for the above-mentioned distances; iii) quantifying the almost sure convergence; iv) relaxing the i.i.d.\! assumption on the $\xitil_i$'s. See, e.g., the books \cite{dudleyCLT,shorack,vdVaartWelln} for a comprehensive treatment of both the original problem and its developments. In particular, the very fruitful line of research---also pursued in the present paper---which aims at providing \emph{quantitative rates} of the Glivenko-Cantelli convergence never ceased to be investigated from the end of the Sixties. See, for example, \cite{BobkovLedoux,boissard,boissardlegouic,dudleyGC,fourn,gozlan1,horowitz,massart,talagrand,weedbach} for an approach based on the empirical processes theory, and \cite{ajtai,ambrosio,caracciolo,dereich,shor,yukich} for the perspective of the so-called optimal matching problem. Possible applications to various areas of pure and applied mathematics can be found in the references of the quoted papers, with a particular mention to: i) statistics (see \cite{dudleyGC} and also \cite{dudleyCLT,efronBOOK,shorack,vdVaartWelln}); ii) particle physics and PDE's (see \cite{ambrosio,bolley,caracciolo,fourn}); iii) numerical analysis, with particular attention to quantization (see \cite{dereich,shor,yukich,weedbach}); iv) machine learning (see also \cite{lugosi,VC}).

The present paper tackles the above quantification problem for various forms of $\hat{\pfrak}_n$ by studying rates with respect to both $L^p$ and almost sure convergence. In Section \ref{sect:iid}, $\pfrak_0$ is fixed and the $\xitil_i$'s are i.i.d.\! from $\pfrak_0$, while in Section \ref{sect:exchangeable} the $\xitil_i$'s are assumed to be exchangeable. In both sections, as in many other papers, the discrepancy between two p.m.'s is measured in terms of the distance $\ud_{[\Space]}^{(p)}$, usually called $p$-\emph{Wasserstein} distance (see however the Bibliographical Notes at the end of Chapter 6 of \cite{villani} for a correct attribution), defined as follows: given a separable metric space $(\Space, \ud_{\Space})$, $p \geq 1$, and $\mu_1, \mu_2$ in $[\Space]_p := \big{\{}\mu\ \mbox{p.m. on}\ (\Space, \mathscr{B}(\Space))\ \big{|}\ \int_{\Space} [\ud_{\Space}(x, x_0)]^p \mu(\ud x) < +\infty\ \mbox{for some}\ x_0 \in \Space\big{\}}$, put
\begin{equation}\label{eq:Wp}
\ud_{[\Space]}^{(p)}(\mu_1, \mu_2) := \inf_{\gamma \in \mathcal{F}(\mu_1, \mu_2)}\Big(
\int_{\Space^2} [\ud_{\Space}(x, y)]^p\ \gamma(\ud x \ud y)\Big)^{1/p}
\end{equation}
where $\mathcal{F}(\mu_1, \mu_2)$ stands for the class of all p.m.'s on $(\Space^2, \mathscr{B}(\Space^2))$ with $i$-th marginal equal to $\mu_i$, $i = 1, 2$. Here and throughout, given any topological space $\mathbb{X}$, 
$\mathscr{B}(\mathbb{X})$ stands for the Borel $\sigma$-algebra on $\mathbb{X}$, and the product $\mathbb{X}^n$ is thought of as endowed with the product topology, for any $n \in \naturals \cup \{\infty\}$. The couple $(\pms_p, \ud_{\pms}^{(p)})$ proves to be a separable metric space, which is complete if $(\Space, \ud_{\Space})$ is also complete, provided that $(\Space, \ud_{\Space})$ satisfies the additional \emph{Radon property}, i.e. every p.m.\! on $(\Space, \mathscr{B}(\Space))$ has the compact inner approximation property. See Definition 5.1.4 and Proposition 7.1.5 in \cite{amgisa}. Finally, recall that a metric subspace $(\Space, \ud_{\Space})$ of another complete and separable metric space, say $(\hat{\Space}, \ud_{\hat{\Space}})$, has the Radon property if $\Space \in \mathscr{B}(\hat{\Space})$. See Theorem 7.1.4 in \cite{du}. 


Throughout the paper, one puts $(\Omega, \mathscr{F}) := (\Space^{\infty}, \mathscr{B}(\Space^{\infty}))$ and defines the sequence $\{\xitil_i\}_{i \geq 1}$ to be the coordinate random elements of $\Space^{\infty}$, i.e.\! $\xitil_i(\omega) := \omega_i$ for all $\omega = (\omega_1,\omega_2, \dots) \in \Space^{\infty}$, providing a formal definition of the sample evoked at the beginning of the section. In addition, one confines oneself to considering forms of $\hat{\pfrak}_n$ which, like $\empiric$, are presentable as measurable functions of $\xinsu$. More formally, given some $p \geq 1$ and defined $\mathscr{F}_n$ to be the sub-$\sigma$-algebra of  $\mathscr{F}$ generated by $\xinsu$, 
the mapping $\hat{\pfrak}_n : (\Omega, \mathscr{F}_n) \rightarrow ([\Space]_p, \mathscr{B}([\Space]_p))$ is required to be measurable for every $n \in \naturals$.

In Section \ref{sect:iid}, one fixes $p \geq 1$, $\pfrak_0$ in $[\Space]_p$ and endows $(\Omega, \mathscr{F})$ with the law $\pfrak_0^{\infty}$ which makes the coordinates i.i.d.\! random elements with common p.d. \!
$\pfrak_0$, i.e.\! $\pfrak_0^{\infty}[\xitil_1 \in A_1, \dots, \xitil_n \in A_n] = \prod_{i=1}^n \pfrak_0^{\infty}[\xitil_i \in A_i] = \prod_{i=1}^n \pfrak_0(A_i)$, for every $n \in \naturals$ and $A_1, \dots, A_n \in \mathscr{B}(\Space)$.
The main objective is to determine an increasing sequence $b_n$ of positive numbers, with $\lim_{n \rightarrow +\infty} b_n = +\infty$, and positive constants $C_p(\pfrak_0)$ and $Y_p(\pfrak_0)$ for which
\begin{equation} \label{eq:Main1}
\pfrak_0^{\infty}\Big[\Big(\sup_{n \geq 1} b_n \ud_{[\Space]}^{(p)}(\pfrak_0, \hat{\pfrak}_n)\Big)^p \Big] \leq C_p(\pfrak_0)
\end{equation}
and
\begin{equation} \label{eq:Main2}
\pfrak_0^{\infty}\Big(\Big\{\limsup_{n \rightarrow +\infty} b_n \ud_{[\Space]}^{(p)}(\pfrak_0, \hat{\pfrak}_n) \leq Y_p(\pfrak_0)\Big\}\Big) = 1\ .
\end{equation}
In (\ref{eq:Main1}), as well as throughout the paper, expectation of a real random variable $X$ with respect to a p.m.\! $\mu$ on $(\Omega, \mathscr{F})$ is denoted by $\mu(X)$. Moreover, in the sequel, the expression ``$\pfrak_0$ satisfies (\ref{eq:Main1})'' (``$\pfrak_0$ satisfies (\ref{eq:Main2})'', respectively) will be used as a shorthand to mean that, given $\pfrak_0$, (\ref{eq:Main1}) ((\ref{eq:Main2}), respectively) is in force with a suitable constant $C_p(\pfrak_0)$ ($Y_p(\pfrak_0)$, respectively). 
\begin{remark} \label{rmk:Wps}
It is worth mentioning how to extend the validity of (\ref{eq:Main1}) to lower exponents $s \in [1,p]$, once (\ref{eq:Main1}) has been established. In fact, one has
$$
\pfrak_0^{\infty}\Big[ \Big( \sup_{n \geq 1} b_n \ud_{[\Space]}^{(s)}(\pfrak_0, \hat{\pfrak}_n) \Big)^s \Big] \leq \Big\{ \pfrak_0^{\infty} \Big[ \Big( \sup_{n \geq 1} b_n \ud_{[\Space]}^{(s)}(\pfrak_0, \hat{\pfrak}_n)\Big)^p \Big] \Big\}^{\frac{s}{p}} 
\!\! \leq \Big\{ \pfrak_0^{\infty} \Big[ \Big( \sup_{n \geq 1} b_n \ud_{[\Space]}^{(p)}(\pfrak_0, \hat{\pfrak}_n) \Big)^p \Big] \Big\}^{\frac{s}{p}}
$$
where the former inequality follows from the Lyapunov inequality for moments, while the latter comes from the monotonicity of the Wasserstein distance with respect to the order. Once (\ref{eq:Main2}) has been established, 
the same monotonicity entails $\pfrak_0^{\infty}\big(\big\{\limsup b_n \ud_{[\Space]}^{(s)}(\pfrak_0, \hat{\pfrak}_n) \leq Y_p(\pfrak_0)\big\}\big) = 1$ for any $s \in [1,p]$.  
\end{remark}

The first result in Section \ref{sect:iid} (Theorem \ref{thm:fournier}) states the validity of (\ref{eq:Main1})-(\ref{eq:Main2}) when $\hat{\pfrak}_n = \empiric$, $\Space = \rd$,
$\ud_{\Space}$ coincides with the standard Euclidean distance, $p \in [1,+\infty) \cap (d/2, +\infty)$ and $b_n \sim (n/\log\log n)^{1/2p}$, provided that $\int_{\rd} |\xb|^{2p+\delta} \pfrak_0(\ud \xb) < +\infty$ obtains for some $\delta > 0$. To introduce the second topic of the same section, one should notice that $b_n$ increases slower and slower as $d$ gets large. As noted in \cite{ajtai,fourn,massart,shor,yukich}, this drawback is an intrinsic feature of the problem since $\lim_{n \rightarrow +\infty} n^{1/d} \pfrak_0^{\infty}\big[\ud_{[\rd]}^{(1)}(\pfrak_0, \empiric)\big] > 0$ when $\pfrak_0$ coincides with the uniform measure on the cube $[-1, 1]^d$. This phenomenon is ascribable to the fact that the use of $\empiric$ in approximating $\pfrak_0$ is thoroughly justified when there is no significant \emph{a priori} restriction on the subclass of $[\Space]_p$ in which $\pfrak_0$ can be fixed, so that the slowdown of the divergence of $b_n$ for large dimensions seems to be the price to pay for this generality. In point of fact, many problems in applied mathematics and statistics provide enough information to restrict the aforesaid admissible class to a family $\model$ of distinguished p.m.'s $\mu_{\theta}$ on $(\Space, \mathscr{B}(\Space))$, determined up to some finite-dimensional parameter $\theta$ in $\Theta \subseteq \rk$ so that $\theta \mapsto \mu_{\theta}$ is injective. In this framework, the approximation of $\pfrak_0 = \mu_{\theta_0}$ proves to be more natural within the elements of $\model$, so that one first approximates $\theta_0$ by a suitable random element $\hat{\theta}_n \in \Theta$, depending on $\xinsu$ according to well-known principles of \emph{statistical estimation}, and then approximates $\pfrak_0$ by $\hat{\pfrak}_n = \mu_{\hat{\theta}_n}$. The second main result supports this last method of approximation, in contrast to the rougher one based on the choice of $\empiric$, by showing that, if $\hat{\pfrak}_n = \mu_{\hat{\theta}_n}$, one can put $p = 2$ and $b_n \sim (n/\log\log n)^{1/2}$ in (\ref{eq:Main1})-(\ref{eq:Main2}) independently of the dimension of the $\xitil_i$'s and of $\Theta$, at least when $\model$ coincides with the class of all non-singular multidimensional Gaussian distributions (Theorem \ref{thm:Gauss}) and, more in general, with a distinguished type of statistical exponential family (Theorem \ref{thm:ExpFam}).

Before explaining Section \ref{sect:exchangeable}, it is worth commenting on the value of the problems associated with (\ref{eq:Main1})-(\ref{eq:Main2}) and their solutions contained in Theorems \ref{thm:fournier}-\ref{thm:ExpFam}. First, the use of the distance $\ud_{[\Space]}^{(p)}$ connects the present problems with the so-called Vapnik-Chervonenkis theory. See \cite{dudleyCLT,vdVaartWelln,VC}. However, the actual novelty of the present study consists in finding rates according to a concept of \emph{strong} (i.e. \emph{uniform} with respect to $n$) convergence, obtained by suitable applications of classical inequalities concerning the \emph{dominated ergodic theorem}, due to Siegmund and Teicher. In fact, the current literature has investigated, until now, only non uniform bounds like $\pfrak_0^{\infty}\Big[\Big(\ud_{[\Space]}^{(p)}(\pfrak_0, \hat{\pfrak}_n)\Big)^p \Big] \leq C^{'}(\pfrak_0)\alpha(n)$, valid for all $n \in \naturals$ and suitable positive constants $C^{'}(\pfrak_0)$, whilst (\ref{eq:Main1}) entails, of course, $\pfrak_0^{\infty}\Big[\Big(\sup_{n \geq N} \ud_{[\Space]}^{(p)}(\pfrak_0, \hat{\pfrak}_n)\Big)^p \Big] \leq C_p(\pfrak_0)/b_N^p$, for all $N \in \naturals$. In spite of the strengthening expressed by the latter inequality, the rate $\alpha(n)$ given in the literature turns out to be comparable with our $1/b_n^p$, in the sense that $\alpha(n)$ goes to zero slightly faster than $1/b_n^p$ under the same assumption $\int_{\rd} |\xb|^{2p+\delta} \pfrak_0(\ud \xb) < +\infty$, at least for small $\delta > 0$. Cf. Theorem 1 in \cite{fourn} and Theorem 13 in \cite{gozlan1}. Moreover, apart from the obvious formal improvement, uniform bounds, which are strictly connected with the concept of \emph{uniform integrability}, prove to be of crucial importance in those applications where one is interested in showing that the union of a collection of ``bad events''---and not only the ``bad events''
taken individually---has small probability. Indeed, one can notice that (\ref{eq:Main1}) entails $\pfrak_0^{\infty}\left(\bigcup_{n \geq 1} \{\ud_{[\Space]}^{(p)}(\pfrak_0, \hat{\pfrak}_n) > M/b_n\} \right) = \pfrak_0^{\infty}\left(\sup_{n \geq 1} b_n \ud_{[\Space]}^{(p)}(\pfrak_0, \hat{\pfrak}_n) > M \right) \leq C_p(\pfrak_0)/M^p$ by means of Markov's inequality, so that, given any $\eta > 0$, one can choose $M_{\eta}$ large enough to obtain $\pfrak_0^{\infty}\big(\cap_{n \geq 1} \{\ud_{[\Space]}^{(p)}(\pfrak_0, \hat{\pfrak}_n) \leq M_{\eta}/b_n\} \big) \geq 1-\eta$. On the other hand, certain applications to mathematical statistics and physics would prove completely justified only in the presence of such uniform bounds, as suggested, for example, by Wiener in relation to the problem of connecting von Neumann's ergodic theorem with Birkhoff's. See Section 1 of \cite{yosida}. As far as our personal motivations to study the uniform convergence displayed in (\ref{eq:Main1}), we mention the paper \cite{cifdorega}, which indicates a line of research---pursued in the next Section \ref{sect:exchangeable}---focused on some approximations of posterior and predictive distributions in Bayesian statistics. Indeed, we realized that it is just a condition of the same type as (\ref{eq:Main1}) that allows the application of a suitable martingale argument due to Blackwell and Dubins, so that we seize the opportunity to highlight this connection.

In Section \ref{sect:exchangeable}, the coordinate random elements $\xitil_1, \xitil_2, \dots$ are assumed to be exchangeable, i.e. the joint distribution of every finite subset of $n$ of these elements depends only on $n$ and not on the particular subset, for every $n \in \naturals$. This condition is obviously satisfied if the $\xitil_i$'s are i.i.d.. Actually, exchangeability is the suitable assumption to be made to translate the usual empirical condition of analogy (not necessarily identity) of the observations into precise probabilistic terms. It is well-known that, given any exchangeable p.m.\! $\rho$ on $(\Omega, \mathscr{F}) = (\Space^{\infty}, \mathscr{B}(\Space^{\infty}))$, there is an extension of the Glivenko-Cantelli theorem, due to de Finetti \cite{definetti30,definettiLincei,definetti37}, which states the existence of a \emph{random} p.m., say $\randommeasure : (\Omega, \mathscr{F}) \rightarrow ([\Space],  \mathscr{B}([\Space]))$, such that $\rho(\{\empiric \Rightarrow \randommeasure,\ \text{as}\ n \rightarrow \infty\}) =1$, $\Rightarrow$ denoting the limit in the sense of weak convergence of p.m.'s, and $[\Space]$ the class of all p.m.'s on $(\Space, \mathscr{B}(\Space))$ endowed with the topology generated by $\Rightarrow$. Moreover, the $\xitil_i$'s turn out to be conditionally i.i.d. given $\randommeasure$, with common p.d. equal to $\randommeasure$. These two facts are encapsulated in the well-known \emph{de Finetti representation theorem}, equivalently reformulated as follows (see \cite{ald}): for any exchangeable p.d.\! $\rho$ on $(\Omega, \mathscr{F})$, there exists one and only one p.m. $\pi$ on $(\pms, \mathscr{B}(\pms))$, equal to the p.d. of $\randommeasure$, such that $\rho(A) = \int_{[\Space]} \pfrak^{\infty}(A) \pi(\ud\pfrak)$ holds for any $A \in \mathscr{F}$, $\pfrak^{\infty}$ denoting the p.m.\! on $(\Omega, \mathscr{F})$ that makes the coordinate random variables i.i.d.\! with common distribution $\pfrak$. In Bayesian statistics, $\pi$ is usually called \emph{prior distribution}. One encounters therein two new questions concerning \emph{posterior} and \emph{predictive} distributions, i.e.\! regular conditional p.d.'s $\pi(\xinsu) := \pi(\xinsu, \cdot)$ and $p(\xinsu) := p(\xinsu, \cdot)$ of $\randommeasure$ and $\{\xitil_{n+i}\}_{i \geq 1}$, respectively, given $\xinsu$. In particular, $\pi(\xinsu, B)$ is a version of $\rho(\randommeasure \in B\ |\ \xinsu)$ and $p(\xinsu, A)$ a vesion of $\rho(\{\xi_{i+n}\}_{i \geq 1} \in A\ |\ \xinsu)$, for $B \in \mathscr{B}(\pms)$ and $A \in \mathscr{F}$. Apropos of these questions, Section \ref{sect:exchangeable} completes and enriches the previous work \cite{cifdorega} by extending some of its statements and by providing a Bayesian interpretation of the main results (\ref{eq:Main1})-(\ref{eq:Main2}) formulated in Section \ref{sect:iid}. As to the approximation of $\pi(\xinsu)$, the aim is to show that, if $\rho[C_p(\randommeasure)] < +\infty$ with the same $C_p$ as in (\ref{eq:Main1}), then
\begin{equation} \label{eq:Main3}
\rho\Big(\Big\{ \limsup_{n \rightarrow +\infty} b_n \ud_{[[\Space]_p]}^{(p)}\big{(} \pi(\xinsu), \delta_{\hat{\pfrak}_n} \big{)} \leq Y_p(\randommeasure) \Big\}\Big) = 1
\end{equation}
holds with the same $Y_p$ as in (\ref{eq:Main2}), where $\ud_{[[\Space]_p]}^{(p)}$ is defined according to (\ref{eq:Wp}) with the proviso that $(\Space, \ud_{\Space})$ therein is replaced by $([\Space]_p, \ud_{[\Space]}^{(p)})$. See Theorem \ref{thm:BD} in Section \ref{sect:exchangeable}. Notice that $\ud_{[[\Space]_p]}^{(p)}\big{(} \pi(\xinsu), \delta_{\hat{\pfrak}_n} \big{)} = \left(\int_{\pms} [\ud_{[\Space]}^{(p)}(\pfrak, \hat{\pfrak}_n)]^p \pi(\xinsu, \ud \pfrak)\right)^{1/p}$. It is also worth remarking that the choice of $\ud_{[[\Space]_p]}^{(p)}$ as reference distance establishes an interesting link between (\ref{eq:Main2}) and (\ref{eq:Main3}), which shows how the latter can be viewed as an extension of the former: when the $\xitil_i$'s are i.i.d.\! with common p.d.\! $\pfrak_0$, $\rho$ coincides with $\pfrak_0^{\infty}$, $\pi(\xinsu)$ is equal to $\delta_{\pfrak_0}$ with $\rho$-probability 1, and $\ud_{[[\Space]_p]}^{(p)}\big{(} \pi(\xinsu), \delta_{\hat{\pfrak}_n} \big{)} = \ud_{\pms}^{(p)}\big{(} \pfrak_0, \hat{\pfrak}_n \big{)}$. The approximation of $p(\xinsu)$ is a more delicate problem, here solved, consistently with the general lines and notation given in \cite{cifdorega,cifdoregaNOTE}, by considering conditional p.d.'s, say $q_m(\xinsu) := q_m(\xinsu, \cdot)$, of $\empiricm := \frac{1}{m} \sum_{i=1}^m \delta_{\xitil_{i+n}}$ given $\xinsu$, for every $m \in \naturals$. In particular, $q_m(\xinsu, B)$ is a version of $\rho(\empiricm \in B\ |\ \xinsu)$ for $B \in \mathscr{B}(\pms)$. Therefore, the main result concerning the approximation of the predictive distributions reads
\begin{equation} \label{eq:Main4}
\rho\Big(\Big\{\limsup_{n \rightarrow +\infty} b_n \ud_{[[\Space]_p]}^{(p)}\big{(} q_m(\xinsu), \hat{\pfrak}_n^{\infty} \circ \empiricm^{-1} \big{)} \leq Y_p(\randommeasure)\Big\}\Big) = 1
\end{equation}
for every $m \in \naturals$, while $\hat{\pfrak}_n^{\infty} \circ \empiricm^{-1}(B)$ stands for $\hat{\pfrak}_n^{\infty}\big(\big\{\frac{1}{m} \sum_{i=1}^m \delta_{\xitil_{i+n}} \in B\big\}\big)$, for $B \in \mathscr{B}(\pms)$. See Theorem \ref{thm:BDpred} in Section \ref{sect:exchangeable}. After establishing these general facts, the next new results, i.e.\! Theorems from \ref{thm:BayesNP} to \ref{thm:BayesExpFam}, reconsider Theorems from \ref{thm:fournier} to \ref{thm:ExpFam}, respectively, according to the present Bayesian perspective. Thus, explicit expressions for $\hat{\pfrak}_n$, $b_n$, $C_p(\randommeasure)$ and $Y_p(\randommeasure)$ are provided in each of the Theorems from \ref{thm:BayesNP} to \ref{thm:BayesExpFam}.

\begin{remark} \label{rmk:BayesWps}
In the same vein of Remark \ref{rmk:Wps}, once (\ref{eq:Main3}) has been established, the identities
$$
\rho\Big(\Big\{ \limsup_{n \rightarrow +\infty} b_n \ud_{[[\Space]_s]}^{(s)}\big{(} \pi(\xinsu), \delta_{\hat{\pfrak}_n} \big{)} \leq Y_p(\randommeasure) \Big\}\Big) = \rho\Big(\Big\{ \limsup_{n \rightarrow +\infty} b_n \ud_{[[\Space]_p]}^{(p)}\big{(} \pi(\xinsu), \delta_{\hat{\pfrak}_n} \big{)} \leq Y_p(\randommeasure) \Big\}\Big) = 1
$$
are valid for any $s \in [1,p]$, since $\ud_{[[\Space]_s]}^{(s)}\big{(} \pi(\xinsu), \delta_{\hat{\pfrak}_n} \big{)} \leq \ud_{[[\Space]_p]}^{(s)}\big{(} \pi(\xinsu), \delta_{\hat{\pfrak}_n} \big{)} \leq \ud_{[[\Space]_p]}^{(p)}\big{(} \pi(\xinsu), \delta_{\hat{\pfrak}_n} \big{)}$. In addition, for any $s \in [1,p]$, (\ref{eq:Main4}) entails $\rho\big(\big\{ \limsup b_n \ud_{[[\Space]_s]}^{(s)}\big{(} q_m(\xinsu), \hat{\pfrak}_n^{\infty} \circ \empiricm^{-1} \big{)} \leq Y_p(\randommeasure)\big\}\big) = 1$.
\end{remark}

The achievement of (\ref{eq:Main3})-(\ref{eq:Main4}) represents a specific asymptotic analysis of both posterior and predictive distributions, which is relevant with a view to a better understanding of different kinds of empirical approximation to orthodox Bayesian methods. In fact, the impressive growth of Bayesian statistics in the last decades has produced new complex models---in particular of nonparametric type---which, although appreciated for their predictive features, generate very often serious hurdles to clear from the point of view of direct computations. These difficulties are very often circumvented by the use of unorthodox, but more manageable, Bayesian methods, such as: empirical Bayes (see \cite{efronBOOK,maritz,morris,robbins}); partial and profile likelihood (see \cite{cox,severini}); numerical techniques based on simulation of random quantities, like the bootstrap (see \cite{efron,efronBOOK}). In particular, contrasting $\hat{\pfrak}_n^{\infty} \circ \empiricm^{-1}$ to $q_m(\xinsu)$ as in (\ref{eq:Main4}) is a cornerstone in bootstrap techniques. Very recently, the question of merging of orthodox and empirical Bayes procedures has been studied in 
\cite{rizzelli,petrone,rousseau}, but with a significant difference: (\ref{eq:Main3})-(\ref{eq:Main4}) are reformulated therein by replacing $\rho$ with some hypothetical distribution $\pfrak_{\star}^{\infty}$, in the spirit of the \emph{what if method} described, e.g., in \cite{diafree1}. Other results which are more comparable to ours, though confined to the merging of predictive distributions with $m=1$, can be found in \cite{BCPR,BPR}.

In conclusion, the ultimate aim of producing bounds like (\ref{eq:Main3})-(\ref{eq:Main4}) is to quantify the degree of accuracy in the approximation ensuing from one of the aforesaid empirical techniques. In fact, for any $\epsilon, \eta > 0$, one can find $n_0 = n_0(\epsilon, \eta)$ and $L = L(\eta)$ such that $\rho\Big[\cap_{n \geq n_0} \{\ud_{[[\Space]_p]}^{(p)}\big{(} \pi(\xinsu), \delta_{\hat{\pfrak}_n} \big{)} \leq (L+\epsilon)/b_n\} \Big] \geq 1-\eta$ and
$\rho\Big[\cap_{n \geq n_0} \{\ud_{[[\Space]_p]}^{(p)}\big{(} q_m(\xinsu), \hat{\pfrak}_n^{\infty} \circ \empiricm^{-1} \big{)} \leq (L+\epsilon)/b_n\} \Big] \geq 1-\eta$ are in force. Thus, one can observe that the determination of $L$ follows from the specific form of $Y_p(\randommeasure)$, while that of $n_0$ represents an interesting open problem to be tackled in future works.

\section{Results for the i.i.d. case} \label{sect:iid}

This section contains four propositions. The first one, of a general character, is a re-statement of an inequality for normed sums of i.i.d. real random variables, originally due to Siegmund \cite{siegmund} and Teicher \cite{teicher}, sometimes referred to as \emph{dominated ergodic theorem}. See also Section 10.3 in \cite{chte}. In the following version, this inequality is reinforced, with respect to the original statements, by the explicit characterization of the upper bound, which plays a crucial role in the proofs of the remaining theorems exhibited in the present section. These last results, in turn, provide affirmative answers to the achievement of (\ref{eq:Main1})-(\ref{eq:Main2}) in three situations of interest.

Hence, one starts with the reformulation of the Siegmund-Teicher inequality.
\begin{prop}\label{prop:teicher}
\emph{Let} $\{X_i\}_{i \geq 1}$ \emph{be a sequence of i.i.d. real random variables on} $\probabilityspace$ \emph{with} $\ee[X_1] = 0$ \emph{and} $\ee[|X_1|^r] < +\infty$ \emph{for some} $r > 2$. \emph{Then, after putting} $\mathrm{Log}_2(x) := \ind\{x < e^e\} + (\log\log x) \ind\{x \geq e^e\}$, \emph{one has}
\begin{equation} \label{eq:teicherpiu2}
\ee\Big[\sup_{n \geq 1} \frac{\big{|}\sum_{i=1}^n X_i\big{|}^r}{(n \mathrm{Log}_2(n))^{r/2}}\Big] \leq \sigma^r \Big[
\alpha_0(r) + \alpha_1(r) \left(\frac{\ee[|X_1|^r]}{\sigma^r}\right)^{\lceil r \rceil}\Big]
\end{equation}
\emph{with suitable constants} $\alpha_0(r), \alpha_1(r)$, $\sigma^2 := \ee[X_1^2]$ \emph{and} $\lceil r \rceil := \inf\{n \in \naturals\ |\ n \geq r\}$.
\end{prop}
\noindent See Subsection \ref{sect:proofteicher} for the proof.

\begin{remark} \label{rmk:teicher}
With a view to successive applications of this proposition, it is crucial to underline that the constants $\alpha_0(r)$ and $\alpha_1(r)$ appearing in the upper bound do not depend on the law of the $X_i$'s. Precise expressions for these constants can be drawn from the combination of specific passages of the proof, culminating in inequality (\ref{eq:teicherFinal}).
\end{remark}

Combination of (\ref{eq:teicherpiu2}) with certain inequalities concerning $\ud_{[\Space]}^{(p)}$, originally proved in \cite{dereich} and reformulated more recently in \cite{fourn}, yields the following
\begin{thm} \label{thm:fournier}
\emph{If} $\Space = \rd$ \emph{and} $\ud_{\Space}$ \emph{coincides with the Euclidean distance, assume that, for some} $\delta > 0$ \emph{and} $p \in [1,+\infty) \cap (d/2, +\infty)$, $\int_{\rd} |\xb|^{2p+\delta} \pfrak_0(\ud \xb) < +\infty$ \emph{obtains. Then}, (\ref{eq:Main1})-(\ref{eq:Main2}) \emph{hold true with} $\hat{\pfrak}_n = \empiric$, $b_n = (n/\mathrm{Log}_2 n)^{1/2p}$. \emph{Moreover, for any choice of} $r$ \emph{in} $(2, 3/M)$, \emph{with} 
$M := \max\{2 - p/d; 1 + 1/(2+\delta/p)\}$, \emph{one can put} 
\begin{eqnarray}
C_p(\pfrak_0) &=& \frac{K(p,d) \overline{C}(r)}{1 - 2^{-\lambda}} \Big[1 + \frac{2^p}{1 - 2^{-\sigma}} \Big(\int_{\rd} |\xb|^{2p+\delta} \pfrak_0(\ud \xb)\Big)^{\frac{3}{r} - 1}\Big] \label{eq:CpNonParametric} \\
Y_p(\pfrak_0) &=& \Big\{\frac{\sqrt{2} K(p,d)}{1 - 2^{-(p - d/2)}} \Big[1 + \frac{2^p}{1 - 2^{-\delta/2}} \Big(\int_{\rd} |\xb|^{2p+\delta} \pfrak_0(\ud \xb)\Big)^{1/2}\Big]\Big\}^{1/p} \label{eq:YpNonParametric}
\end{eqnarray}
\emph{where} $K(p,d), \overline{C}(r), \lambda$ \emph{and} $\sigma$ \emph{are suitable positive constants independent of} $\pfrak_0$.
\end{thm}
\noindent See Subsection \ref{sect:ProofFournier} for the proof, as well as for the exact evaluation of $K(p,d), \overline{C}(r), \lambda$ and $\sigma$. It is also worth recalling Remark \ref{rmk:Wps} in order to indicate how to enlarge the applicability of this theorem, since one could be interested in bounding the Wasserstein distance of low order (e.g., 1 or 2) even in the presence of a high dimension $d$.  

Before introducing the other results, it is worth noticing once again the troublesome dependence, in the previous theorem, of $b_n$ and $p$ on $d$, which appears to be unavoidable when $\pfrak_0$ is viewed as any element of the whole class $[\rd]_{2p+\delta}$. See \cite{ajtai,fourn,massart,shor,yukich,weedbach}. Then, one is naturally led to search for situations in which the aforesaid drawback does not occur. As already remarked in Section \ref{sect:intro}, there are many significant problems in which it is natural to constrain both $\pfrak_0$ and $\hat{\pfrak}_n$ within a distinguished class $\model$ strictly contained in $[\Space]$. By the way of illustration, maintaining the choice $\Space = \rd$, one now considers the noteworthy case in which $\model = \{\mu_{\theta}\}_{\theta \in \Theta}$ coincides with the class of the non-degenerate Gaussian p.d.'s, namely
\begin{equation}\label{eq:Gauss}
\mu_{\theta}(\ud \xb) = \left(\frac{1}{2\pi}\right)^{d/2} (\text{det}(V))^{-1/2} \exp\left\{-\frac{1}{2}\ ^t(\xb - \mb)V^{-1}(\xb - \mb)\right\}\ud \xb\ \ \ \ \ (\xb \in \rd)
\end{equation}
where $\theta = (\mb, V)$, $\Theta = \rd \times \mathbb{GL}_{sym}^{+}(d)$ and $\mathbb{GL}_{sym}^{+}(d)$ denotes the class of symmetric and positive-definite $d \times d$ matrices with real entries. As for the approximating p.d. $\hat{\pfrak}_n$, in agreement with the so-called \emph{plug-in} method of statistical point estimation, put $\hat{\pfrak}_n = \mu_{\hat{\theta}_n}$,
$\hat{\theta}_n = (\hat{\mb}_n, \hat{V}_n)$ and
\begin{equation}\label{eq:MLEGauss}
\hat{\mb}_n = \frac{1}{n} \sum_{i=1}^{n} \xitil_i, \ \ \ \ \ \hat{V}_n = \frac{1}{n} \sum_{i=1}^{n} (\xitil_i - \hat{\mb}_n)\ ^t(\xitil_i - \hat{\mb}_n)\ ,
\end{equation}
with the proviso that $\mu_{\hat{\theta}_1} = \delta_{\xitil_1}$, to obtain
\begin{thm} \label{thm:Gauss}
\emph{Let} $\Space = \rd$ \emph{and} $\ud_{\Space}$ \emph{coincide with the Euclidean distance. Assume that} $\pfrak_0 = \mu_{\theta_0}$, $\mu_{\theta_0}$ \emph{being defined by} (\ref{eq:Gauss}) \emph{with} $\theta_0 = (\mb_0, V_0) \in \rd \times \mathbb{GL}_{sym}^{+}(d)$. \emph{Then,} $\pfrak_0$ \emph{satisfies} (\ref{eq:Main1})-(\ref{eq:Main2}) \emph{with} $p = 2$, $b_n = (n/\mathrm{Log}_2 n)^{1/2}$, $\hat{\pfrak}_n = \mu_{\hat{\theta}_n}$ \emph{and}
$\hat{\theta}_n = (\hat{\mb}_n, \hat{V}_n)$ \emph{as in} (\ref{eq:MLEGauss}). \emph{Moreover, for any} $\varepsilon$ \emph{such that} $0 < \varepsilon \leq [2(d+1)]^{-d}$, \emph{there exist two positive constants} $K(\varepsilon,d)$ \emph{and} $H(\varepsilon,d)$ \emph{for which one can put}
\begin{eqnarray}
C_2(\pfrak_0) &=& [c_{\star} + b_d^2(d+1)] \text{tr}[V_0] + c_{\ast}K(\varepsilon,d) \|V_0\|_F + H(\varepsilon,d) \label{eq:CpGauss} \\
Y_2(\pfrak_0) &=& \sqrt{2} \big\{\sigma_{max}[V_0] + d\sqrt{K(\varepsilon,d)} \sqrt{\text{tr}[V_0]} \big\} , \label{eq:YpGauss}
\end{eqnarray}
\emph{where} $\text{tr}[\cdot]$, $\|\cdot\|_F$ \emph{and} $\sigma_{max}^2[V_0]$ \emph{stand for the trace, the Frobenius norm and the maximum eigenvalue of a matrix, respectively, while} $c_{\star}$ \emph{and} $c_{\ast}$ \emph{are two numerical constants, even independent of the dimension} $d$. 
\end{thm}
\noindent See Subsection \ref{sect:proofGauss} for the proof. In particular, the proof of Lemma \ref{lm:kullback} contains the precise evaluation of $K(\varepsilon,d)$, while $c_{\star}$ is displayed in (\ref{eq:TeicherMeanGauss}), $c_{\ast}$ in (\ref{eq:lastGaussInt}), and $H(\varepsilon,d)$ in (\ref{eq:lastGauss}). Lastly, definition and properties of $\|\cdot\|_F$ can be found in Section 5.6 of \cite{horn}. \\

A remarkably interesting fact is that an analogous conclusion holds true for a rather general type of \emph{exponential family} which, besides including the previous class of Gaussian distributions as a distinguished, but particular, case, enjoys important, global properties from the point of view of mathematical statistics. Two of them are worth noticing: i) identifiability, in the sense of a parametric family of p.m.'s; ii) existence and uniqueness of the maximum likelihood estimator, expressed as arithmetic mean of values of a given function $\tb$. In turn, i)-ii) are implied by certain weak conditions of a technical character, which are explained, for example, in \cite{barndorff,brown}. A brief summary is given here, referring the reader to the quoted monographs for the details. One can start introducing a $\sigma$-finite measure $\mu$ on the metric space $(\Space, \ud_{\Space})$, to dominate the whole family $\model$ (to be defined just below), together with the aforesaid measurable map $\tb : (\Space, \mathscr{B}(\Space)) \rightarrow (\rk, \mathscr{B}(\rk))$, in such a way that:
\begin{enumerate}
\item[A)] the interior $\Theta$ of the convex hull of the support of $\mu \circ \tb^{-1}$ is non-void,
\item[B)] the set $\Lambda := \Big{\{}\yb \in \rk\ \Big{|}\ \int_{\Space} e^{\yb \cdot \tb(x)} \mu(\ud x) < +\infty\Big{\}}$ is a non void open subset of $\rk$, which proves to be convex.
\end{enumerate}
Thus, a parametrization of the exponential family $\mathcal{M}$, conducive to i)-ii), is given by
\begin{equation}\label{eq:meanvalueparametrization}
\mu _{\thetab}(\ud x) := \exp\{V^{-1}(\thetab) \cdot \tb(x) - M(V^{-1}(\thetab))\} \mu(\ud x)\ \ \ \ \ (\thetab \in \Theta)
\end{equation}
where $M(\yb) := \log\left(\int_{\Space} e^{\yb \cdot \tb(x)} \mu(\ud x)\right) : \rk \rightarrow (-\infty, +\infty]$ and $V(\yb) := \nabla M(\yb) = \int_{\Space} \tb(x) \exp\{\yb \cdot \tb(x) - M(\yb)\} \mu(\ud x)$. These definitions make sense in view of the following facts: first, Theorems 1.13, 2.2 and 2.7 in \cite{brown} state that $M$ is strictly convex on $\Lambda$, lower semi-continuous on $\rk$, of class $C^{\infty}(\Lambda)$ and analytic. Second, Corollary 5.3 in \cite{barndorff} implies that $M$ is steep (essentially smooth), according to Definition 3.2 of \cite{brown}.  Third, from Theorem 3.6 in \cite{brown}, $\yb \mapsto V(\yb)$ defines a smooth homeomorphism of $\Lambda$ and $\Theta$. Thus, putting 
 $\gamma_{\yb}(A) := \int_A  \exp\{\yb \cdot \tb(x) - M(\yb)\} \mu(\ud x)$, $A \in \mathscr{B}(\Space)$, Corollary 2.5 in \cite{brown} entails the equivalence of $\gamma_{\yb_1} = \gamma_{\yb_2}$, $\yb_1 = \yb_2$ and $\int \tb \ud \gamma_{\yb_1} = \int \tb \ud \gamma_{\yb_2}$, making precise the requisite of identifiability. Therefore, (\ref{eq:meanvalueparametrization}) defines the family $\model$ which is the subject of the next theorem, in such a way that
$\pfrak_0$ will be a fixed element $\mu_{\thetab_0}$ of $\model$, and $\hat{\pfrak}_n$ the element $\mu_{\hat{\thetab}_n}$ corresponding to
\begin{equation}  \label{eq:MLE}
\hat{\thetab}_n := \frac{1}{n} \sum_{i=1}^{n} \tb(\xitil_i) \in \Theta
\end{equation}
which represents, in this case, the \emph{maximum likelihood estimator} of the parameter $\thetab$.

The particular form of $\hat{\thetab}_n$ is conducive to an application of Proposition \ref{prop:teicher}, but with the significant variant that the analog of the supremum appearing in (\ref{eq:teicherpiu2}) is multiplied by a function, say $\Phi$, of $\hat{\thetab}_n$, which diverges as either $|\hat{\thetab}_n|$ goes to infinity or $\hat{\thetab}_n$ approaches the boundary of $\Theta$. Hence, to obtain any quantitative bound, the following concepts will come in usefull:
first, the cumulant generating function, say $\Psi$, of $\tb(\xitil_1)$, namely $\Psi(\yb) := \log\left(\int_{\Space} e^{\yb \cdot \tb(x)} \mu_{\thetab_0}(\ud x)\right)$. Second, the \emph{Legendre transformation}, say $I_{\mu}$, of $M$, that is
\begin{equation}  \label{eq:Imu}
I_{\mu}(\thetab) := \sup_{\yb \in \Lambda} \{\thetab \cdot \yb - M(\yb)\} = V^{-1}(\thetab) \cdot \thetab - M(V^{-1}(\thetab))\ .
\end{equation}
Third, the Legendre transformation, say $I_{\thetab_0}$, of $\Psi$, given by
\begin{equation}  \label{eq:Itheta0}
I_{\thetab_0}(\thetab) := \sup_{\yb \in \Lambda} \{\thetab \cdot \yb - \Psi(\yb)\} = I_{\mu}(\thetab) - \thetab \cdot \yb_0 - M(\yb_0)
\end{equation}
where $\yb_0 = V^{-1}(\thetab_0)$. See Section 26 of \cite{rockafellar} (especially Theorems 26.4 and 26.6) for well-definiteness and regularity of $I_{\mu}$, as well as Theorem 6.9 in \cite{brown} for steepness of $I_{\mu}$. The way is now paved to specify $\Phi$, i.e.
\begin{equation} \label{eq:PHI}
\Phi(\etab) := \sqrt{\int_0^1 (1-s) \|\textrm{Hess}[I_{\mu}](\thetab_0 + s(\etab - \thetab_0))\|_F\ \ud s}\ \ \ \ \ \ (\etab \in \Theta)
\end{equation}
where $\|\cdot\|_F$ stands for the Frobenius norm. As a consequence of the properties of $I_{\mu}$, in order that $\Phi(\etab)$ might diverge it is necessary that either $|\etab|$ goes to infinity or $\etab$ approaches 
$\partial\Theta$. Hence, with a view to quantifying the rapidity of the divergence of $\Phi$, one observes that there always exist two functions $\Phi_1 : [0, +\infty) \rightarrow [0, +\infty)$ and $\Phi_2 : \Theta \rightarrow [0, +\infty)$ satisfying
\begin{equation} \label{eq:PHI1-2}
\left\{
\begin{array}{l}
\Phi(\etab) \leq \Phi_1(|\etab|) + \Phi_2(\etab)\ \text{for\ all\ } \etab \in \Theta; \\
\Phi_1\ \text{is\ increasing\ and\ diverges\ at\ infinity\ }; \\
\forall\ \{\etab_n\}_{n \geq 1} \subset \Theta, \Phi_2(\etab_n) \rightarrow +\infty\ \text{entails\ } \ud(\etab_n, \partial\Theta) \rightarrow 0,\ \text{as}\ n \rightarrow +\infty\ .
\end{array} \right.
\end{equation}
Then, along with $\Phi_2$, one introduces the non-negative integer (possibly equal to $+\infty$) $N_{\Phi_2}(\rho, \sigma)$ which represents the minimum number of convex closed sets that form a covering of $R_1 := \{\thetab \in \Theta\ |\ \Phi_2(\thetab) \geq \frac{\sigma}{2}, |\thetab - \thetab_0| \leq \rho\}$ with the constrain that such a covering must be included in $R_2 := \{\thetab \in \Theta\ |\ \Phi_2(\thetab) \geq \frac{\sigma}{4}, |\thetab - \thetab_0| \leq \rho\}$.

Finally, it should be noticed that the involvement of $\ud_{[\Space]}^{(p)}$ in (\ref{eq:Main1})-(\ref{eq:Main2}) forces to assume that $\model$ is included in $[\Space]_p$ and, with a view to the obtainment of explicit upper bounds, induces to focus attention on the so-called \emph{Talagrand inequality}. With reference to the elements of $\model$ with $p=2$, this inequality reads
\begin{equation} \label{eq:Talagrand}
\ud_{[\Space]}^{(2)}(\mu_{\thetab_0}, \mu_{\thetab})^2 \leq C_T(\thetab_0) \mathrm{K}(\mu_{\thetab}\ |\ \mu_{\thetab_0}) \ \ \ \ \ ( \thetab_0, \thetab \in \Theta)
\end{equation}
where $\mathrm{K}(\mu_{\thetab}\ |\ \mu_{\thetab_0}) := \int_{\Space} \log\left(\frac{\ud \mu_{\thetab}}{\ud \mu_{\thetab_0}}(x)\right)\mu_{\thetab}(\ud x)$ is the \emph{Kullback-Leibler information} of $\mu_{\thetab_0}$ at $\mu_{\thetab}$. See, e.g., \cite{gozlanAOP,gozlanSurvey} for further information.

\begin{thm} \label{thm:ExpFam}
\emph{Suppose that} A)-B) \emph{are fulfilled and that the elements of the exponential family} (\ref{eq:meanvalueparametrization})
\emph{belong to} $[\Space]_2$ \emph{and meet} (\ref{eq:Talagrand}). \emph{Then,} $\pfrak_0 = \mu_{\thetab_0}$ \emph{satisfies} (\ref{eq:Main2}) \emph{with} $p=2$,  $b_n = (n/\mathrm{Log}_2 n)^{1/2}$, $\hat{\pfrak}_n = \mu_{\hat{\thetab}_n}$, $\hat{\thetab}_n$ \emph{as in} (\ref{eq:MLE}), \emph{and}
\begin{equation} \label{eq:Yexp}
Y_2(\pfrak_0) = \sqrt{2 C_T(\thetab_0)} \Phi(\thetab_0) \sigma_{max}(\thetab_0)
\end{equation}
\emph{where} $\sigma_{max}^2(\thetab_0)$ \emph{is the largest eigenvalue of} $\mathrm{Hess}[M](V^{-1}(\thetab_0))$. \emph{Moreover,}
$\pfrak_0$ \emph{satisfies} (\ref{eq:Main1}) \emph{with the same} $p$, $b_n$, $\hat{\pfrak}_n$ \emph{and a suitable} $C_2(\pfrak_0)$ \emph{specified in the proof if, in addition, there exist a positive constant} $\tau(\thetab_0)$ \emph{and functions} $\rho, \sigma : [\tau(\thetab_0), +\infty) \rightarrow [0,+\infty)$ \emph{satisfying}
\begin{equation} \label{eq:CostExp}
\int_{\tau(\thetab_0)}^{+\infty} m(t)^{-2} e^{-m(t)} \ud t + \int_{\tau(\thetab_0)}^{+\infty} N_{\Phi_2}(\rho(t), \sigma(t)) h(t)^{-2} e^{-h(t)} \ud t < +\infty
\end{equation}
\emph{where, with reference to} (\ref{eq:PHI1-2}), \emph{one puts} $m(t) := \min\{\rho(t), \frac{t}{\sigma(t)}, \Phi_1^{-1}(\frac{\sigma(t)}{2})\}$ \emph{and} $h(t) := \inf\{I_{\thetab_0}(\thetab)\ |\\ \Phi_2(\thetab) \geq \frac{\sigma(t)}{4}, |\thetab - \thetab_0| \leq \rho(t)\}$.
\end{thm}
\noindent  See Subsection \ref{sect:proofExpFam} for the proof, which also includes the elements to specify the constant $C_2(\pfrak_0)$.

Before concluding the section, it is worth remarking that the abundance of assumptions in the last theorem counterbalances the generality of the family (\ref{eq:meanvalueparametrization}). In any case, checking of A)-B) is a standard task, as shown in Chapter 1 of \cite{brown}. To check the validity of (\ref{eq:Talagrand}) see \cite{gozlanAOP,gozlanSurvey}, which also highlight the interesting connection with the \emph{logarithmic Sobolev inequality} as means to estimate the constant $C_T(\thetab_0)$. Lastly, the finding of $\tau(\thetab_0)$, $\rho(t)$, $\sigma(t)$ and $N(t)$ is less standard and could prove more labored, even because of lack of background literature.

\section{Exchangeable random variables} \label{sect:exchangeable}

As announced in Section \ref{sect:intro}, the random coordinates $\xitil_1, \xitil_2, \dots$ are now thought of as exchangeable random elements distributed according to the p.m.\! $\rho$ on $(\Omega, \mathscr{F}) = (\Space^{\infty}, \mathscr{B}(\Space^{\infty}))$. The reader is recommended to resort to the representation theorem recalled therein and, in particular, to the meaning of the random p.m.\! $\randommeasure$ as limit of the sequence of the empirical laws $\empiric$'s, as well as of the symbols $\pi(\xinsu)$ and $p(\xinsu)$ to denote the posterior and the predictive distribution, respectively.

The first theorem deals with the posterior distribution, establishing the validity of (\ref{eq:Main3}). Its novelty is that it provides information about the \emph{range} of the sequence $\{\ud_{[\pms_p]}^{(p)}(\pi(\xisun), \delta_{\hat{\pfrak}_n})\}_{n \geq 1}$, when $(\Space, \ud_{\Space})$ is a general separable metric space meeting the Radon property. 
\begin{thm} \label{thm:BD}
\emph{Assume that, for some} $p \geq 1$, $\hat{\pfrak}_n \in [\Space]_p$ \emph{for all} $n \in \naturals$ \emph{and} $\randommeasure \in \pms_p$ \emph{with} $\rho$-\emph{probability} 1. \emph{Moreover, let} $\{b_n\}_{n \geq 1}$ \emph{be an increasing and diverging sequence of positive numbers with respect to which} $\randommeasure$ \emph{satisfies} (\ref{eq:Main1})-(\ref{eq:Main2}) \emph{with} $\rho$-\emph{probability} 1. \emph{Then,}
\begin{equation} \label{eq:MainP}
\rho\big[\limsup_{n \rightarrow \infty} b_n \ud_{[\pms_p]}^{(p)}\big(\pi(\xisun), \delta_{\hat{\pfrak}_n}\big) \leq Y_p(\randommeasure) \big] = 1
\end{equation}
\emph{provided that} $C_p(\randommeasure)$ \emph{and} $Y_p(\randommeasure)$ \emph{are real random variables and} $\rho[C_p(\randommeasure)] <+\infty$.
\end{thm}
\noindent The proof is deferred to Subsection \ref{sect:proofBD}. Apropos of the assumption that $C_p(\randommeasure)$ and $Y_p(\randommeasure)$ are real random variables, notice that their checking boils down to an analysis of their explicit expressions. Cf., e.g., (\ref{eq:CpNonParametric})-(\ref{eq:YpNonParametric}) and (\ref{eq:CpGauss})-(\ref{eq:YpGauss}). Lastly, recall Remark \ref{rmk:BayesWps} for an answer to the problem of bounding $\limsup b_n \ud_{[\pms_s]}^{(s)}\big(\pi(\xisun), \delta_{\hat{\pfrak}_n}\big)$ for $s \in [1,p]$.

In the same vein, one now deals with the approximation of the distributions $q_m(\xinsu)$ as in (\ref{eq:Main4}). To connect $q_m(\xinsu)$ with $\pi(\xinsu)$, suffice it to notice that $p(\xinsu, A) = \int_{\pms} \pfrak^{\infty}(A) \pi(\xinsu, \ud \pfrak)$ and $q_m(\xinsu, B) = p(\xinsu, \empiricm^{-1}(B))$ hold true for every $A \in \mathscr{B}(\Space^{\infty})$ and $B \in \mathscr{B}(\pms)$. From a technical point of view, knowledge of the latter distribution may help to make the properties of the former more explicit. This is what happens when, by following the usual frequentistic interpretation of analogy of the observations, one approximates $q_m(\xinsu)$ by means of the distinguished distribution of $\frac{1}{m} \sum_{i=1}^m \delta_{\xitil_{i+n}}$ when the $\xitil_{i+n}$'s are viewed as i.i.d.\! with common p.d.\! $\hat{\pfrak}_n$, for every $n$. Apropos of this, one states the second main result of this section as
\begin{thm} \label{thm:BDpred}
\emph{If, for some} $p \geq 1$, $\hat{\pfrak}_n \in [\Space]_p$ \emph{for all} $n \in \naturals$ \emph{and} $\randommeasure \in \pms_p$ \emph{with} $\rho$-\emph{probability} 1, \emph{then}
\begin{equation} \label{eq:savare1}
\rho\left[\ud_{[[\Space]_p]}^{(p)}(q_m(\xinsu), \hat{\pfrak}_n^{\infty} \circ \empiricm^{-1}) \leq \ud_{[[\Space]_p]}^{(p)}(\pi(\xinsu), \delta_{\hat{\pfrak}_n})\right] = 1
\end{equation}
\emph{for every} $m \in \naturals$. \emph{In particular, under the same assumptions of} Theorem \ref{thm:BD}, \emph{one has}
\begin{equation} \label{eq:MainPred}
\rho\left[\limsup_{n \rightarrow \infty} b_n \ud_{[[\Space]_p]}^{(p)}(q_m(\xinsu), \hat{\pfrak}_n^{\infty} \circ \empiricm^{-1}) \leq Y_p(\randommeasure)\right] = 1
\end{equation}
\emph{for every} $m \in \naturals$.
\end{thm}
\noindent See Subsection \ref{sect:proofBDpred} for the proof, and resort to Remark \ref{rmk:BayesWps} for an answer to the problem of bounding Wasserstein distances of lower orders $s \in [1,p]$.

Theorems \ref{thm:BD} and \ref{thm:BDpred} are deliberately stated in a very general form, as far as the nature of the space $\Space$ is concerned, but contain the strong assumption that (\ref{eq:Main1})-(\ref{eq:Main2}) are in force. With a view to their real use, one could do a good turn by providing possibly simple sufficient conditions for the validity of  (\ref{eq:Main1})-(\ref{eq:Main2}), even if one has to give up the generality of $\Space$. It is with this spirit that we state the following three theorems which, from a merely mathematical side, are immediate consequences of the combination of Theorems \ref{thm:BD} and \ref{thm:BDpred} with Theorems \ref{thm:fournier}, \ref{thm:Gauss} and \ref{thm:ExpFam}, respectively. The first statement deals with the nonparametric model, and can be viewed as an extension of Theorem 4 in \cite{cifdorega} from $\reals$ to $\rd$.
\begin{thm} \label{thm:BayesNP}
\emph{Let} $\Space = \rd$ \emph{and} $\ud_{\Space}$ \emph{be the Euclidean distance. Assume that, for some} $\delta > 0$ \emph{and} $p \in [1,+\infty) \cap (d/2, +\infty)$, $\rho\Big{[} \Big{(} \int_{\rd} |\xb|^{2p+\delta} \randommeasure(\ud\xb) \Big{)}^{\frac{3}{r} - 1} \Big{]} = \int_{[\rd]} \big{(} \int_{\rd} |\xb|^{2p+\delta} \pfrak(\ud\xb) \big{)}^{\frac{3}{r} - 1} \pi(\ud\pfrak) < +\infty$
\emph{is valid with any suitable} $r \in (2,3)$ \emph{for which} $p/d > 2 - \frac{3}{r}$ \emph{and} $(2p+\delta)(\frac{3}{r} - 1) > p$. \emph{Then}, (\ref{eq:Main3})-(\ref{eq:Main4}) \emph{hold true with} $\hat{\pfrak}_n = \empiric$, $b_n = (n/\mathrm{Log}_2 n)^{1/2p}$ \emph{and the same constant} $Y_p$ \emph{as in} Theorem \ref{thm:fournier}, \emph{with} $\randommeasure$ \emph{in place of} $\pfrak_0$.
\end{thm}

The next theorem provides a new result concerning the $d$-dimensional Gaussian model.
\begin{thm} \label{thm:BayesGauss}
\emph{If} $\Space = \rd$ \emph{and} $\ud_{\Space}$ \emph{coincides with the Euclidean distance, suppose that} $\randommeasure = \mu_{\tilde{\theta}}$ \emph{with} $\rho$-\emph{probability} 1, \emph{where} $\tilde{\theta} = (\tilde{\mathbf{m}}, \tilde{V})$ \emph{is a random element taking values in} $\Theta = \rd \times \mathbb{GL}^{+}_{sym}(d)$ \emph{and} $\mu_{\theta}$ \emph{is referred to} (\ref{eq:Gauss}). \emph{If} 
$\rho\big[\|\tilde{V}\|_F\big] = \int_0^{+\infty} x \pi_{\|\tilde{V}\|_F}(\ud x) < +\infty$ \emph{is valid, with} $\pi_{\|\tilde{V}\|_F}$ \emph{standing for the p.d. of} $\|\tilde{V}\|_F$, \emph{then} (\ref{eq:Main3})-(\ref{eq:Main4}) \emph{hold true with} $\hat{\pfrak}_n = \mu_{\hat{\theta}_n}$, $b_n = (n/\mathrm{Log}_2 n)^{1/2}$ \emph{and the same constant} $Y_p$ \emph{as in} Theorem \ref{thm:Gauss}, \emph{with} $\tilde{V}$ \emph{in place of} $V_0$.
\end{thm}

Finally, the third result concerns the exponential family considered in Theorem \ref{thm:ExpFam}.
\begin{thm} \label{thm:BayesExpFam}
\emph{Let} $\randommeasure$ \emph{equal} $\mu_{\tilde{\thetab}}$ \emph{with} $\rho$-\emph{probability} 1, \emph{where} $\mu_{\thetab}$ \emph{is referred to} (\ref{eq:meanvalueparametrization}) \emph{and} $\tilde{\thetab}$ \emph{is a random element with values in} $\Theta$. \emph{Suppose that} A)-B) \emph{are fulfilled and that, with} $\rho$-\emph{probability} 1, $\mu_{\tilde{\thetab}} \in \pms_2$, (\ref{eq:Talagrand}) \emph{holds with} $\tilde{\thetab}$ \emph{in place of} $\thetab_0$ \emph{and} $\hat{\thetab}_n$ \emph{in place of} $\thetab$, \emph{and} (\ref{eq:CostExp}) \emph{is in force. Then, if the same constant} $C_2(\pfrak_0)$ \emph{mentioned in} Theorem \ref{thm:ExpFam} \emph{is such that} $C_2(\mu_{\tilde{\thetab}})$ \emph{is a real random variable satisfying} $\rho[C_2(\mu_{\tilde{\thetab}})] < +\infty$, \emph{properties} (\ref{eq:Main3})-(\ref{eq:Main4}) \emph{hold true with} $\hat{\pfrak}_n = \mu_{\hat{\theta}_n}$, $b_n = (n/\mathrm{Log}_2 n)^{1/2}$ \emph{and} $Y_2(\randommeasure) = \sqrt{2 C_T(\tilde{\thetab})} \Phi(\tilde{\thetab}) \sigma_{max}(\tilde{\thetab})$.
\end{thm}

\begin{remark}
Theorems \ref{thm:BD} and \ref{thm:BDpred} face the challenging task of comparing probability laws of p.d.'s. Nevertheless, there are situations in which it suffices to have a partial knowledge of a p.d., such as mean values or appropriate synthetic measures of interesting aspects of that p.d., presentable as functionals on classes of p.d.'s. Thus, Theorems \ref{thm:BD} and \ref{thm:BDpred} are ideally accompanied by results concerning the comparison between $\rho[g(\randommeasure)\ |\ \xinsu]$ and $g(\hat{\pfrak}_n)$, and between $\rho[g(\frac{1}{m} \sum_{i=1}^m \delta_{\xitil_{i+n}})\ |\ \xinsu]$ and $\hat{\pfrak}_n^{\infty}[g(\frac{1}{m} \sum_{i=1}^m \delta_{\xitil_{i+n}})]$, respectively, where $g$ stands for some specific functional. In particular, in view of the Kantorovich-Rubinstein representation of $\ud^{(1)}$ (see, e.g., Theorem 11.8.2 of \cite{du}), 
if $g$ is Lipschitz-continuous, one has
$$
\big{|} \rho[g(\randommeasure)\ |\ \xinsu] - g(\hat{\pfrak}_n) \big{|} \leq \text{Lip}(g) \ud_{[\pms_1]}^{(1)}\left(\pi(\xisun), \delta_{\hat{\pfrak}_n}\right)
$$
and
$$
\Big{|} \rho\Big[g\Big(\frac{1}{m} \sum_{i=1}^m \delta_{\xitil_{i+n}}\Big)\ |\ \xinsu\Big] - \hat{\pfrak}_n^{\infty}\Big[g\Big(\frac{1}{m} \sum_{i=1}^m \delta_{\xitil_{i+n}}\Big)\Big]  \Big{|} \leq \text{Lip}(g) \ud_{[[\Space]_1]}^{(1)}(q_m(\xinsu), \hat{\pfrak}_n^{\infty} \circ \empiricm^{-1})
$$
where $\text{Lip}(g)$ stands for the Lipschitz constant of $g$. Due to its relevance, proven by a rich literature, the empirical approximation of functionals will be studied in a forthcoming paper.
\end{remark}

\section{Proofs} \label{sect:proofs}

Gathered here are the proofs of the main results.

\subsection{Proof of Proposition \ref{prop:teicher}} \label{sect:proofteicher}

This subsection complements the arguments already developed to prove Lemma 1 and Theorem 4 in Section 10.3 of \cite{chte}, in order to justify the bound (\ref{eq:teicherpiu2}). Then, notation is substantially the same as therein, and numbering relates to that very same book.

To start, consider a sequence $\{Y_n\}_{n \geq 1}$ of independent, non-negative random variables, and set $\mathfrak{S} := \sum_{n=1}^{+\infty} Y_n$ and $\mathfrak{M}_k := \sum_{n=1}^{+\infty} \ee[Y_n^k]$, for $k > 0$. An application of the inequality $(a + b)^m \leq 2^{m-1}(a^m + b^m)$, valid for $a, b > 0$ and $m \geq 1$, yields $\ee[\mathfrak{S}^q] \leq 2^{q-2}(\mathfrak{M}_q + \mathfrak{M}_1 \ee[\mathfrak{S}^{q-1}])$ for $q \in \{2, 3, \dots\}$. Thus, by induction, one gets
$\ee[\mathfrak{S}^q] \leq \lambda_{q, q-1} \mathfrak{M}_1^q + \sum_{k=0}^{q-2} \lambda_{q, k} \mathfrak{M}_{q-k} \mathfrak{M}_1^k$
for $q \in \{2, 3, \dots\}$, where $\lambda_{2,1} := 1$, $\lambda_{q, 0} := 2^{q-2}$ and $\lambda_{q+1, h} := 2^{q-1} \lambda_{q, h-1}$ if $h = 1, \dots, q$. If $q > 2$ and $q \not\in \naturals$, one can write
\begin{eqnarray}
\ee[\mathfrak{S}^q] &\leq& 2^{[q]-1}\left(\mathfrak{M}_q + \mathfrak{M}_{q-[q]} \ee[\mathfrak{S}^{[q]}]\right) \nonumber \\
&\leq& 2^{[q]-1}\Big{(}\mathfrak{M}_q + \lambda_{[q], [q]-1} \mathfrak{M}_{q-[q]} \mathfrak{M}_1^{[q]} + \mathfrak{M}_{q-[q]} \sum_{k=0}^{[q]-2} \lambda_{[q], k} \mathfrak{M}_{[q]-k} \mathfrak{M}_1^k\Big{)}\ . \nonumber
\end{eqnarray}

Now, the reasoning continues parallelling the method used to prove Theorem 4. First, assume that the distribution of the $X_n$'s is symmetric with $\ee[X_1^2] = 1$. Second, define $c_n := (n \mathrm{Log}_2(n))^{-1/2}$, $X_n^{'} := X_n \ind\{|X_n| \leq n^{1/r}\}$,
$X_n^{''} := X_n \ind\{|X_n| > n^{1/r}\}$ and decompose $S_n := \sum_{i=1}^n X_i$ as $S_n = S_n^{'} + S_n^{''}$, with $S_n^{'} := \sum_{i=1}^n X_i^{'}$ and $S_n^{''} := \sum_{i=1}^n X_i^{''}$. Apropos of $S_n^{''}$, exploit the monotonicity of the sequence $\{c_n\}_{n \geq 1}$, as in the last inequality on page 389, to obtain
$$
\ee\Big{[} \sup_{n \geq 1} \frac{|S_n^{''}|^r}{(n \mathrm{Log}_2(n))^{r/2}} \Big{]} \leq \ee\Big{[}\Big{(} \sup_{n \geq 1} c_n \sum_{i=1}^n |X_i^{''}| \Big{)}^r\Big{]} \leq \ee\Big{[}\Big{(} \sum_{n \geq 1} c_n |X_n^{''}| \Big{)}^r\Big{]} \ .
$$
Therefore, if $r \in \{3, 4, \dots\}$, one gets
\begin{eqnarray}
&& \ee\Big{[} \sup_{n \geq 1} \frac{|S_n^{''}|^r}{(n \mathrm{Log}_2(n))^{r/2}} \Big{]} \nonumber \\
&\leq& \lambda_{r, r-1} \left(\sum_{n=1}^{+\infty} c_n \ee[|X_n^{''}|]\right)^r + \sum_{k=0}^{r-2} \lambda_{r, k} \left(\sum_{n=1}^{+\infty} c_n^{r-k} \ee[|X_n^{''}|^{r-k}]\right) \cdot \left(\sum_{n=1}^{+\infty} c_n \ee[|X_n^{''}|]\right)^k
\label{eq:Theorem4Teichera}
\end{eqnarray}
whilst, if $r > 2$ and $r \not\in \naturals$, the following bound
\begin{eqnarray}
&& \ee\Big{[} \sup_{n \geq 1} \frac{|S_n^{''}|^r}{(n \mathrm{Log}_2(n))^{r/2}} \Big{]} \nonumber \\
&\leq& 2^{[r]-1}\Big{\{}\sum_{n=1}^{+\infty} c_n^r \ee[|X_n^{''}|^r] + \lambda_{[r], [r]-1}
\Big{(}\sum_{n=1}^{+\infty} c_n^{r-[r]} \ee[|X_n^{''}|^{r-[r]}]\Big{)} \cdot \Big{(}\sum_{n=1}^{+\infty} c_n \ee[|X_n^{''}|]\Big{)}^{[r]} \nonumber \\
&+& \Big{(}\sum_{n=1}^{+\infty} c_n^{r-[r]} \ee[|X_n^{''}|^{r-[r]}]\Big{)} \cdot \sum_{k=0}^{[r]-2} \lambda_{[r], k}
\Big{(}\sum_{n=1}^{+\infty} c_n^{[r] - k} \ee[|X_n^{''}|^{[r]-k}]\Big{)} \cdot
\Big{(}\sum_{n=1}^{+\infty} c_n \ee[|X_n^{''}|]\Big{)}^k \Big{\}} \ \ \ \ \ \label{eq:Theorem4Teicherb}
\end{eqnarray}
is in force. Now, combining the inequalities of H\"{o}lder and Markov, one gets $\ee[|X_n^{''}|^h] \leq (1/n)^{\frac{r-h}{r}} \ee[|X_1|^r]$ for every $n \in \naturals$ and $h \in (0, r]$ which, in turn, yields $\sum_{n=1}^{+\infty} c_n^h \ee[|X_n^{''}|^h] \leq K(h) \ee[|X_1|^r]$, with $K(h) := \sum_{n=1}^{+\infty} (1/n)^{\frac{r-h}{r}} c_n^h < +\infty$. Thus, one concludes that
\begin{equation} \label{eq:teicherEsterno}
\ee\Big{[} \sup_{n \geq 1} \frac{|S_n^{''}|^r}{(n \mathrm{Log}_2(n))^{r/2}} \Big{]} \leq \beta_1(r) (\ee[|X_1|^r])^{\lceil r \rceil}
\end{equation}
is valid with a suitable numerical constant $\beta_1(r)$, determined as follows. First, recall that $\ee[X_1^2] = 1$ entails $\ee[|X_1|^r] \geq 1$. Then, if $r \in \{3, 4, \dots\}$, invoke (\ref{eq:Theorem4Teichera}) to get $\beta_1(r) = \sum_{k=0}^{r-1} \lambda_{r, k} K(r-k) (K(1))^k$, whilst, if $r > 2$ and $r \not\in \naturals$, use (\ref{eq:Theorem4Teicherb}) to obtain
$$
\beta_1(r) = 2^{[r]-1}\Big{\{} K(r) + \lambda_{[r], [r]-1} K(r - [r]) (K(1))^{[r]} + K(r - [r])\sum_{k=0}^{[r]-2} \lambda_{[r], k}
K([r] - k) (K(1))^k \Big{\}} \ .
$$

Apropos of the term involving $S_n^{'}$, one starts by observing that
\begin{equation} \label{eq:choquet}
\ee\left[\left(\sup_{n \geq 1} c_n |S_n^{'}| \right)^r\right] \leq u_0^r + r\lambda^r \int_{u_0/\lambda}^{+\infty} u^{r-1} \pp\left[\sup_{n \geq 1} c_n |S_n^{'}| > \lambda u\right] \ud u
\end{equation}
holds for every $\lambda, u_0 > 0$. Then, after putting $n_k := [e^k]$ for $k \in \naturals_0$, formula (15) on page 390 can be invoked to write
\begin{equation} \label{eq:levy}
\pp\left[\sup_{n \geq 1} c_n |S_n^{'}| > \lambda u\right] \leq 4 \sum_{k=0}^{+\infty} \pp\left[c_{n_k} S_{n_{k+1}}^{'} >
\lambda u\right]\ .
\end{equation}
Moreover, at the same page, it is proved that $\pp[S_n^{'} > x] \leq \exp\{-tx + nt^2\}$ is valid for every $n \in \naturals$, $x > 0$ and $t \in [0, n^{-1/r}]$. Then, setting $n = n_{k+1}$, $\gamma_r := \sup_{k \in \naturals} \left(\frac{\mathrm{Log}_2(n_{k+1})}{n_{k+1}}\right)^{1/2} (n_{k+1})^{1/r} <+\infty$, $t = \gamma_r^{-1} \left(\frac{\mathrm{Log}_2(n_{k+1})}{n_{k+1}}\right)^{1/2}$ and $x = \frac{u}{\gamma_r c_{n_k}}$ yields
\begin{equation} \label{eq:teicher33}
\pp\left[c_{n_k} S_{n_{k+1}}^{'} > \frac{u}{\gamma_r} \right] \leq \exp\left\{-\frac{(\alpha u - 1)}{\gamma_r^2} \mathrm{Log}_2(n_{k+1})\right\}
\end{equation}
where $\alpha := \inf_{k \in \naturals_0} \left(\frac{n_k \mathrm{Log}_2(n_{k})}{n_{k+1} \mathrm{Log}_2(n_{k+1})}\right)^{1/2} > 0$. Putting $\lambda = \gamma_r^{-1}$, the combination of (\ref{eq:choquet}), (\ref{eq:levy}) and (\ref{eq:teicher33}) leads to
$$
\ee\left[\left(\sup_{n \geq 1} c_n |S_n^{'}| \right)^r\right] \leq u_0^r + 4r\gamma_r^{-r} \sum_{k=0}^{+\infty} \ \int_{u_0 \gamma_r}^{+\infty} u^{r-1} \exp\left\{-\frac{(\alpha u - 1)}{\gamma_r^2}\mathrm{Log}_2(n_{k+1})\right\} \ud u\ .
$$
To conclude, it is enough to choose $u_0$ large enough so that $\alpha \gamma_r u_0 > 4 \gamma_r^2 + 1$ is in force. In fact, this choice entails, for every $u \geq u_0 \gamma_r$, $\alpha u > 4 \gamma_r^2 + 1$ and $\frac{(\alpha u - 1)}{\gamma_r^2}\mathrm{Log}_2(n_{k+1}) \geq 2 \mathrm{Log}_2(n_{k+1}) + \frac{2\alpha}{4\gamma_r^2 + 1} u$ which, in turn, give
\begin{equation} \label{eq:teicherInterno}
\ee\left[\left(\sup_{n \geq 1} c_n |S_n^{'}| \right)^r\right] \leq u_0^r + 4r\gamma_r^{-r}\left(\sum_{k=0}^{+\infty} \exp\{-2 \mathrm{Log}_2(n_{k+1})\}\right) \cdot \left(\frac{4\gamma_r^2 + 1}{2\alpha}\right)^r \Gamma(r) =: \beta_0(r)\ .
\end{equation}

The proof can be now carried out by means of the following steps. First, combine (\ref{eq:teicherEsterno}) with (\ref{eq:teicherInterno}) to get
\begin{eqnarray}
\ee\left[\left(\sup_{n \geq 1} c_n |S_n| \right)^r\right] &\leq& 2^{r-1} \ee\left[\left(\sup_{n \geq 1} c_n |S_n^{'}| \right)^r\right] +
2^{r-1} \ee\left[\left(\sup_{n \geq 1} c_n |S_n^{''}| \right)^r\right] \nonumber \\
&\leq& 2^{r-1}\left(\beta_0(r) + \beta_1(r)(\ee[|X_1|^r])^{\lceil r \rceil}\right)\ , \label{eq:teicherFinal}
\end{eqnarray}
which is valid for symmetric $X_n$'s with $\ee[X_1^2] = 1$. Then, invoke the symmetrization argument at the end of page 390 to show that this last very same bound holds also for non-symmetric $X_n$'s, provided that $\ee[X_1] = 0$ and $\ee[X_1^2] = 1$. Finally, remove the assumption $\ee[X_1^2] = 1$, by replacing every $X_i$ with $X_i/\sigma$.

\subsection{Proof of Theorem \ref{thm:fournier}} \label{sect:ProofFournier}

Fix $d \in \naturals$, $p \in [1, +\infty) \cap (d/2, +\infty)$ and $\delta > 0$. After setting $\beta := 2p + \delta$, observe that it is possible to choose $r \in (2,3)$ so that $\lambda = \lambda(p,d,r) := p - d(2 - \frac{3}{r})$ and $\sigma = \sigma(p,d,r) := \beta(\frac{3}{r} - 1) - p$ are both strictly positive constants. In particular, the bound on $\sigma$ is equivalent to $r < \frac{3\beta}{\beta+p} \in (2,3)$. As for the bound on $\lambda$, notice that it is trivially satisfied for any choice of $r \in (2,3)$ if $p \geq d$ whilst, if $p \in [1, +\infty) \cap (d/2, d)$, one gets $r < \frac{3d}{2d - p} \in (2,3)$. Therefore, these arguments show that $r$ can be chosen as any element of $(2, 3/M)$, where $M := \max\{2 - p/d; 1 + 1/(2+\delta/p)\}$.

Now, consider the LHS of (\ref{eq:Main1}) with $b_n^p = \sqrt{n/\mathrm{Log}_2 n}$ and invoke Lemmata 5 and 6 in \cite{fourn}, whose main contents are now briefly recalled. For $l \in \naturals_0 := \{0, 1, 2, \dots\}$, denote by $\mathcal{P}_l$ the natural partition of $B_0 := (-1, 1]^d$ into $2^{dl}$ translations of $(-2^{-l}, 2^{-l}]^d$ and, for $s \in \naturals$, put $B_s := (-2^s, 2^s]^d\setminus (-2^{s-1}, 2^{s-1}]^d$ and, for any $F \subseteq \rd$, $2^s F := \{2^s \xb\ |\ \xb \in F\}$. In this notation, for any $p \geq 1$ and $\mu,\nu \in [\rd]_p$,
\begin{equation} \label{eq:Fournier}
\ud_{[\rd]}^{(p)}(\mu, \nu)^p \leq K(p,d) \sum_{s=0}^{+\infty} 2^{ps} \sum_{l=0}^{+\infty} 2^{-lp} \sum_{F \in \mathcal{P}_l} \big{|} \mu(2^s F \cap B_s) - \nu(2^s F \cap B_s) \big{|}
\end{equation}
holds with $K(p,d) := 2^{2p} d^{p/2} (2^p + 1)(2^p - 1)^{-2}$. Whence,
\begin{gather}
\pfrak_0^{\infty}\Big[ \sup_{n \geq 1} b_n^p \ud_{[\rd]}^{(p)}(\pfrak_0, \empiric)^p \Big] \nonumber \\
\leq K(p,d)\sum_{s=0}^{+\infty} 2^{ps} \sum_{l=0}^{+\infty} 2^{-lp} \sum_{F \in \mathcal{P}_l} \pfrak_0^{\infty} \Big[\sup_{n \geq 1} b_n^p \Big{|} \pfrak_0(2^s F \cap B_s) - \frac{1}{n}\sum_{i=1}^n \ind\{\xitil_i \in 2^s F \cap B_s\} \Big{|}\Big]\ . \nonumber
\end{gather}
The objective is now to bound the term $\pfrak_0^{\infty}\left[\sup_{n \geq 1} b_n^p \Big{|} \pfrak_0(2^s F \cap B_s) - \frac{1}{n}\sum_{i=1}^n \ind\{\xitil_i \in 2^s F \cap B_s\} \Big{|}\right]$ by the following steps: first, choose $r \in (2,3)$ satisfying the restrictions imposed in the first lines of this proof. Second, apply Lyapunov's inequality to get
\begin{gather}
\pfrak_0^{\infty}\Big[ \sup_{n \geq 1} b_n^p \Big{|} \pfrak_0(2^s F \cap B_s) - \frac{1}{n}\sum_{i=1}^n \ind\{\xitil_i \in 2^s F \cap B_s\} \Big{|}\Big] \nonumber \\
\leq \Big\{\pfrak_0^{\infty}\Big[ \Big( \sup_{n \geq 1} b_n^p \Big{|} \pfrak_0(2^s F \cap B_s) - \frac{1}{n}\sum_{i=1}^n \ind\{\xitil_i \in 2^s F \cap B_s\} \Big{|} \Big)^r \Big] \Big\}^{1/r}\ . \nonumber
\end{gather}
At this stage, the identities $\pfrak_0^{\infty}[\ind\{\xitil_i \in 2^s F \cap B_s\}] = \pfrak_0(2^s F \cap B_s)$ and
$$
b_n^p \Big{|} \pfrak_0(2^s F \cap B_s) - \frac{1}{n}\sum_{i=1}^n \ind\{\xitil_i \in 2^s F \cap B_s\} \Big{|} = \frac{\big{|} \sum_{i=1}^n [\ind\{\xitil_i \in 2^s F \cap B_s\} - \pfrak_0(2^s F \cap B_s)]\big{|}}{\sqrt{n \mathrm{Log}_2 n}}\ ,
$$
valid for any $i$ and $n$ in $\naturals$, respectively, pave the way for the application of Proposition \ref{prop:teicher} with $X_i = \ind\{\xitil_i \in 2^s F \cap B_s\} - \pfrak_0(2^s F \cap B_s)$. Notice that $\sigma^2$ and $\ee[|X_1|^r]$ appearing in the RHS of (\ref{eq:teicherpiu2}) are equal to $\pfrak_0(2^s F \cap B_s)(1 - \pfrak_0(2^s F \cap B_s))$ and $\pfrak_0(2^s F \cap B_s)(1 - \pfrak_0(2^s F \cap B_s))^r + \pfrak_0(2^s F \cap B_s)^r(1 - \pfrak_0(2^s F \cap B_s))$, respectively. In view of $(A+B)^{1/r} \leq A^{1/r} + B^{1/r}$, one can bound from above the ($1/r$)-power of the RHS of (\ref{eq:teicherpiu2}) by means of
\begin{eqnarray}
&& \alpha_0(r)^{1/r} \sqrt{\pfrak_0(2^s F \cap B_s)(1 - \pfrak_0(2^s F \cap B_s))} \nonumber \\
&+& \alpha_1(r)^{1/r} [\pfrak_0(2^s F \cap B_s)(1 - \pfrak_0(2^s F \cap B_s))]^{\frac{3}{r} - 1} \cdot [(1 - \pfrak_0(2^s F \cap B_s))^{r-1} + \pfrak_0(2^s F \cap B_s)^{r-1}]^{\frac{3}{r} - 1} \nonumber
\end{eqnarray}
so that, taking account of $(1 - \pfrak_0(2^s F \cap B_s))^{r-1} + \pfrak_0(2^s F \cap B_s)^{r-1} \leq 2$, $1 - \pfrak_0(2^s F \cap B_s) \leq 1$ and $\frac{1}{2} > \frac{3}{r} - 1$, one can write
$$
\pfrak_0^{\infty} \Big[\sup_{n \geq 1} b_n^p \Big{|} \pfrak_0(2^s F \cap B_s) - \frac{1}{n}\sum_{i=1}^n \ind\{\xitil_i \in 2^s F \cap B_s\} \Big{|} \Big] \leq \overline{C}(r) \left[\pfrak_0(2^s F \cap B_s)\right]^{\frac{3}{r} - 1}\ .
$$
with $\overline{C}(r) := \alpha_0(r)^{1/r} + 2^{3/r}\alpha_1(r)^{1/r}$. At this stage, after recalling that $\sharp \mathcal{P}_l = 2^{dl}$ and that $\frac{3}{r} - 1 \in (0,1)$, another application of Lyapunov's inequality yields
\begin{eqnarray}
&& \sum_{F \in \mathcal{P}_l} \pfrak_0^{\infty}\left[\sup_{n \geq 1} b_n^p \Big{|} \pfrak_0(2^s F \cap B_s) - \frac{1}{n}\sum_{i=1}^n \ind\{\xitil_i \in 2^s F \cap B_s\} \Big{|}\right] \nonumber \\
&\leq& \overline{C}(r) 2^{dl(2 - \frac{3}{r})} \left(\sum_{F \in \mathcal{P}_l} \pfrak_0(2^s F \cap B_s) \right)^{\frac{3}{r} - 1}
= \overline{C}(r) 2^{dl(2 - \frac{3}{r})} [\pfrak_0(B_s)]^{\frac{3}{r} - 1}\ . \nonumber
\end{eqnarray}
Whence,
\begin{eqnarray}
&& \sum_{l=0}^{+\infty} 2^{-lp} \sum_{F \in \mathcal{P}_l} \pfrak_0^{\infty}\left[\sup_{n \geq 1} b_n^p \Big{|} \pfrak_0(2^s F \cap B_s) - \frac{1}{n}\sum_{i=1}^n \ind\{\xitil_i \in 2^s F \cap B_s\} \Big{|}\right] \nonumber \\
&\leq& \overline{C}(r) [\pfrak_0(B_s)]^{\frac{3}{r} - 1} \sum_{l=0}^{+\infty} 2^{dl(2 - \frac{3}{r}) - lp} = \frac{\overline{C}(r)}{1 - 2^{-\lambda}} [\pfrak_0(B_s)]^{\frac{3}{r} - 1}\ . \nonumber
\end{eqnarray}
Now, observe that, for any $s \in \naturals$, the definition of $B_s$ entails $\pfrak_0(B_s) \leq \pfrak_0^{\infty}[\{|\xitil_1| > 2^{s-1}\}]$ so that, from Markov's inequality, $\pfrak_0(B_s) \leq 2^{-\beta(s-1)} \pfrak_0^{\infty}[|\xitil_1|^{\beta}]$ holds with the same $\beta$ as defined at the beginning of this proof. To conclude the proof of (\ref{eq:Main1}), gather the above inequalities together to obtain
$$
\pfrak_0^{\infty}\left[\sup_{n \geq 1} b_n^p \ud_{[\rd]}^{(p)}(\pfrak_0, \empiric)^p \right] \leq \frac{K(p,d) \overline{C}(r)}{1 - 2^{-\lambda}} \left[1 + \frac{2^p}{1 - 2^{-\sigma}} \left(\pfrak_0^{\infty}[|\xitil_1|^{\beta}]\right)^{\frac{3}{r} - 1}\right]
$$
which provides the value of $C_p(\pfrak_0)$ displayed in (\ref{eq:CpNonParametric}).

To prove (\ref{eq:Main2}), use again (\ref{eq:Fournier}) to write
\begin{gather}
\limsup_{n \rightarrow +\infty} b_n^p \ud_{[\rd]}^{(p)}(\pfrak_0, \empiric)^p  \nonumber \\
\leq K(p,d)\sum_{s=0}^{+\infty} 2^{ps} \sum_{l=0}^{+\infty} 2^{-lp} \sum_{F \in \mathcal{P}_l} \limsup_{n \rightarrow +\infty} b_n^p \Big{|} \pfrak_0(2^s F \cap B_s) - \frac{1}{n}\sum_{i=1}^n \ind\{\xitil_i \in 2^s F \cap B_s\} \Big{|}\ . \nonumber
\end{gather}
Since, from Hartman-Wintner's law of the iterated logarithm,
$$
\limsup_{n \rightarrow +\infty} b_n^p \Big{|} \pfrak_0(2^s F \cap B_s) - \frac{1}{n}\sum_{i=1}^n \ind\{\xitil_i \in 2^s F \cap B_s\} \Big{|} = \sqrt{2 \pfrak_0(2^s F \cap B_s)(1 - \pfrak_0(2^s F \cap B_s))}
$$
one can argue exactly as above to obtain
$$
\limsup_{n \rightarrow +\infty} b_n^p \ud_{[\rd]}^{(p)}(\pfrak_0, \empiric)^p \leq \frac{\sqrt{2} K(p,d)}{1 - 2^{-(p - d/2)}} \left[1 + \frac{2^p}{1 - 2^{-\delta/2}} \left(\pfrak_0^{\infty}[|\xitil_1|^{\beta}]\right)^{1/2}\right]
$$
which provides the value of $Y_p(\pfrak_0)$ displayed in (\ref{eq:YpNonParametric}).

\subsection{Proof of Theorem \ref{thm:Gauss}} \label{sect:proofGauss}

To facilitate the proof, one premises the following technical result 

\begin{lm} \label{lm:kullback}
\emph{Let} $\lambda_1, \dots, \lambda_d$ \emph{be the eigenvalues (counted with their multiplicity) of a symmetric, positive definite} $d \times d$ \emph{matrix} $M$ \emph{with real entries. If} $\mathrm{det}(M) \geq \varepsilon$ \emph{for some} $\varepsilon \in (0,1)$, \emph{then there exists a positive constant} $K(\varepsilon, d)$ \emph{such that}
\begin{equation} \label{eq:eigenKullback}
\sum_{j=1}^d (\lambda_j - 1 - \log\lambda_j) \leq K(\varepsilon, d) \sum_{j=1}^d (\lambda_j - 1)^2 = K(\varepsilon, d)\|M - \mathrm{Id}_d\|^2_F
\end{equation}
\emph{where} $\mathrm{Id}_d$ \emph{is the} $d\times d$ \emph{identity matrix and} $\|\cdot\|_F$ \emph{stands for the Frobenius norm}. 
\end{lm}

\emph{Proof of the Lemma}. Fix $\delta \in (0,1)$ and define $d_1(M,\delta)$ to be the number of those $\lambda_j$'s belonging to $[1-\delta, 1+\delta]$. Since $K_1(\delta) := \sup_{x \in [1-\delta, 1+\delta]\setminus\{1\}} \frac{(x-1)^2}{x - \log x - 1} > 0$, if $d_1(M,\delta) > 0$ one gets
$$
\sum_{\substack{j \in \{1,\dots,d\} \\ \lambda_j \in [1-\delta, 1+\delta]}} (\lambda_j - 1 - \log \lambda_j) \leq K_1(\delta) \!\!\!\! \sum_{\substack{j \in \{1,\dots,d\} \\ \lambda_j \in [1-\delta, 1+\delta]}} (\lambda_j - 1)^2\ .
$$
If $d_1(M,\delta) = d$, the lemma is proved upon putting $K(\varepsilon, d) = K_1(\delta)$, with $(1-\delta)^d = \varepsilon$. Thus, assume that $d_2 = d_2(M,\delta) := d - d_1(M,\delta) > 0$, and notice that
$\sum_{(\ast)} (\lambda_j - 1 - \log \lambda_j)\ \leq \sum_{(\ast)} \lambda_j - \ \log \prod_{(\ast)} \lambda_j$, where $(\ast)$ is a shorthand to indicate that sums or products run over those $j$'s in $\{1,\dots,d\}$ for which 
$\lambda_j \not\in [1-\delta, 1+\delta]$. Moreover, by assumption, there holds
$$
\prod_{(\ast)} \lambda_j \ = \  \prod_{j=1}^d \lambda_j \ \cdot   \Big(\prod_{\substack{j \in \{1,\dots,d\} \\ \lambda_j \in [1-\delta, 1+\delta]}} \lambda_j \Big)^{-1} \geq \frac{\varepsilon}{(1+\delta)^d} =: p_0
$$
where $p_0 \in (0,1)$. Lyapunov's inequality entails $\frac{1}{d_2} \sum_{(\ast)} \lambda_j \leq \Big( \frac{1}{d_2} \sum_{(\ast)} \lambda_j^2 \ \Big)^{1/2}$ and the GAM inequality gives 
$\frac{1}{d_2} \sum_{(\ast)} \lambda_j \geq \Big( \prod_{(\ast)} \lambda_j \Big)^{1/d_2}$, so that the combination of these inequalities yields 
$p_0^{1/d_2} \leq \Big( \frac{1}{d_2} \sum_{(\ast)} \lambda_j^2 \Big)^{1/2}$. Hence, recalling that $\log p_0 < 0$, there holds
$$
\sum_{(\ast)} (\lambda_j - 1 - \log \lambda_j)\ \leq\ d_2 \Big( \frac{1}{d_2} \sum_{(\ast)} \lambda_j^2 \ \Big)^{1/2} \!\!\! - \log p_0 \
\leq\ \Big[ \Big( \frac{1}{p_0} \Big)^{1/d_2}  -  \frac{1}{d_2} \Big(\frac{1}{p_0}\Big)^{2/d_2} \log p_0 \Big] \ \sum_{(\ast)} \lambda_j^2 \ . \nonumber 
$$
Since $K_2(\delta) := \sup_{x \not\in [1-\delta, 1+\delta]\setminus\{1\}} \frac{x^2}{(x-1)^2} > 0$ for any $\delta > 0$, one deduces the former inequality in (\ref{eq:eigenKullback}). Finally, the latter inequality follows from the properties of the Frobenius norm, as shown, e.g., Section 5.6 of \cite{horn}.\ \ $\square$

The way is now paved for the proof of Theorem \ref{thm:Gauss}. \\

From a well-known expression of the distance $\ud_{[\rd]}^{(2)}$ between two multivariate Gaussian distributions (see, e.g., \cite{dowson,olkin}), one gets
\begin{equation} \label{eq:olkin}
\left[\ud_{[\rd]}^{(2)}(\mu_{\theta_0}, \mu_{\hat{\theta}_n})\right]^2 = |\mathbf{m}_0 - \hat{\mathbf{m}}_n|^2 + \text{tr}\left[V_0 + \hat{V}_n -2(V_0 \hat{V}_n)^{1/2}\right]
\end{equation}
where the symbol $\mu_{\theta}$ is referred to (\ref{eq:Gauss}), while $\hat{\mathbf{m}}_n $ and $\hat{V}_n$ are referred to (\ref{eq:MLEGauss}). Then, to prove (\ref{eq:Main1}), one starts by analyzing the first summand on the above RHS. 
The combination of Lyapunov's inequality with (\ref{eq:teicherpiu2}) yields, for any $r > 2$,
\begin{eqnarray}
\mu_{\theta_0}^{\infty} \Big[ \Big( \sup_{n \geq 1} b_n |\mathbf{m}_0 - \hat{\mathbf{m}}_n | \Big)^2 \Big] &\leq& \sum_{j=1}^d \mu_{\theta_0}^{\infty} \Big[ \Big( \sup_{n \geq 1} \frac{\big{|} \sum_{i=1}^n (\xitil_i^{(j)} - m_0^{(j)}) \big{|}}{\sqrt{n \textrm{Log}_2 n}} \Big)^2 \Big] \nonumber \\
&\leq& \sum_{j=1}^d \Big\{  \mu_{\theta_0}^{\infty} \Big[ \Big( \sup_{n \geq 1} \frac{\big{|} \sum_{i=1}^n (\xitil_i^{(j)} - m_0^{(j)}) \big{|}}{\sqrt{n \textrm{Log}_2 n}} \Big)^r \Big] \Big\}^{2/r} \nonumber \\
&\leq& \sum_{j=1}^d \sigma_j^2 \Big\{ \Big[ \alpha_0(r) + \alpha_1(r) \Big(\frac{\mu_{\theta_0}^{\infty}[|\xitil_1^{(j)} - m_0^{(j)}|^r]}{\sigma_j^r}  \Big)^{\lceil r \rceil} \Big] \Big\}^{2/r} \nonumber \\
&=& \Big[ \alpha_0(r) + \alpha_1(r)\Big( \frac{2^{r/2}}{\sqrt{\pi}} \Gamma\Big(\frac{r+1}{2}\Big) \Big)^{\lceil r \rceil} \Big]^{2/r} \text{tr}[V_0]\ . \ \ \ \label{eq:TeicherMeanGauss}
\end{eqnarray}
As to the latter summand on the RHS of (\ref{eq:olkin}), one can start by writing
$$
\mu_{\theta_0}^{\infty} \Big[ \sup_{n \in \{1, \dots, d\}} b_n^2 \text{tr}[V_0 + \hat{V}_n -2(V_0 \hat{V}_n)^{1/2}] \Big] \leq b_d^2(d+1) \text{tr}[V_0]
$$
since, for every $n \in \naturals$, $\mu_{\theta_0}^{\infty} \Big[\text{tr}[\hat{V}_n]\Big] \leq \text{tr}[V_0]$. Consequently, one can confine oneself to studying $\mu_{\theta_0}^{\infty}\big[ \sup_{n \geq d+1} b_n^2 \text{tr}[V_0 + \hat{V}_n -2(V_0 \hat{V}_n)^{1/2}] \big]$. Now, it is worth observing that the term $\text{tr}\big[V_0 + \hat{V}_n -2(V_0 \hat{V}_n)^{1/2}\big]$ coincides with the squared distance $[\ud_{[\rd]}^{(2)}]^2$ between two $d$-dimensional Gaussian distributions with zero means and a covariance matrix equal to $V_0$ or $\hat{V}_n$, respectively. This fact paves the way for the application of the Talagrand inequality (\ref{eq:Talagrand}), to get
\begin{eqnarray}
\text{tr}\left[V_0 + \hat{V}_n -2(V_0 \hat{V}_n)^{1/2}\right] &\leq&
C_T(\theta_0)\left[ \text{tr}\left( V_0^{-1}\hat{V}_n - \mathrm{Id}_d\right) - \log \text{det}(V_0^{-1}\hat{V}_n) \right] \nonumber \\
&=& C_T(\theta_0) \sum_{j=1}^d [\lambda_j^{(n)} - \log \lambda_j^{(n)} - 1] \label{eq:TalagrandGauss}
\end{eqnarray}
where $\lambda_1^{(n)}, \dots, \lambda_d^{(n)}$ denote the strictly positive (with $\mu_{\theta_0}$-probability 1) eigenvalues of $V_0^{-1}\hat{V}_n$. Note that, as a consequence of the original formulation of the Talagrand inequality (see \cite{talagrandIneq}), $C_T(\theta_0)$ can be chosen equal to $2\|V_0\|_F$. At this stage, after fixing $\varepsilon \in (0,1)$, an application of Lemma \ref{lm:kullback} gives
\begin{gather}
\mu_{\theta_0}^{\infty}\left[\sup_{n \geq d+1} b_n^2 \text{tr}[V_0 + \hat{V}_n -2(V_0 \hat{V}_n)^{1/2}] \ind\{\text{det}(V_0^{-1}\hat{V}_n) \geq \varepsilon\}\right] \nonumber \\
\leq 2 K(\varepsilon, d) \|V_0^{-1}\|_F \cdot \sum_{l,h = 1}^d \mu_{\theta_0}^{\infty}\Big[ \Big(\sup_{n \geq d+1}
\frac{\big{|} \sum_{i=1}^n (\tilde{\eta}_i^{(l,h)} - (V_0)_{l,h}) \big{|}}{\sqrt{n \textrm{Log}_2 n}} \Big)^2 \Big]
\end{gather}
where $\tilde{\eta}_i^{(l,h)} := (\xitil_i^{(l)} - \hat{m}_n^{(l)})(\xitil_i^{(h)} - \hat{m}_n^{(h)})$. Elementary algebra shows that
$\sum_{i=1}^n (\tilde{\eta}_i^{(l,l)} - (V_0)_{l,l}) = \sum_{i=1}^n [(\xitil_i^{(l)} - m_0^{(l)})^2 - (V_0)_{l,l}] - \frac{1}{n}\left(\sum_{i=1}^n [\xitil_i^{(l)} - m_0^{(l)}]\right)^2$ and, for $l \neq h$, $\sum_{i=1}^n (\tilde{\eta}_i^{(l,h)} - (V_0)_{l,h}) = \sum_{i=1}^n [(\xitil_i^{(l)} - m_0^{(l)})(\xitil_i^{(h)} - m_0^{(h)}) - (V_0)_{l,h}] - \frac{1}{n}\left(\sum_{i=1}^n [\xitil_i^{(l)} - m_0^{(l)}] \right) \left(\sum_{i=1}^n [\xitil_i^{(h)} - m_0^{(h)}] \right)$. Whence,
\begin{eqnarray}
\mu_{\theta_0}^{\infty} \Big[ \Big( \sup_{n \geq d+1} \frac{\big{|} \sum_{i=1}^n (\tilde{\eta}_i^{(l,l)} - (V_0)_{l,l}) \big{|}}{\sqrt{n \textrm{Log}_2(n)}} \Big)^2 \Big] &\leq& 2 \mu_{\theta_0}^{\infty} \Big[ \Big( \sup_{n \geq d+1} \frac{\big{|} \sum_{i=1}^n [(\xitil_i^{(l)} - m_0^{(l)})^2 - (V_0)_{l,l}] \big{|}}{\sqrt{n \textrm{Log}_2(n)}} \Big)^2 \Big] \nonumber \\
&+& 2 \mu_{\theta_0}^{\infty} \Big[ \sup_{n \geq d+1} \frac{\Big( \sum_{i=1}^n [\xitil_i^{(l)} - m_0^{(l)}] \Big)^4}{n^3 \textrm{Log}_2(n)} \Big] \label{eq:vardecll}
\end{eqnarray}
and, for $l \neq h$,
\begin{eqnarray}
&& \mu_{\theta_0}^{\infty} \Big[ \Big( \sup_{n \geq d+1} \frac{\big{|} \sum_{i=1}^n (\tilde{\eta}_i^{(l,h)} - (V_0)_{l,h}) \big{|}}{\sqrt{n \textrm{Log}_2(n)}} \Big)^2 \Big] \nonumber \\
&\leq& 2 \mu_{\theta_0}^{\infty} \Big[ \Big( \sup_{n \geq d+1} \frac{\big{|} \sum_{i=1}^n [(\xitil_i^{(l)} - m_0^{(l)})(\xitil_i^{(h)} - m_0^{(h)}) - (V_0)_{l,h}] \big{|}}{\sqrt{n \textrm{Log}_2(n)}} \Big)^2 \Big] \nonumber \\
&+& 2 \mu_{\theta_0}^{\infty} \Big[ \sup_{n \geq d+1} \frac{\Big( \sum_{i=1}^n [\xitil_i^{(l)} - m_0^{(l)}] \Big)^2 \Big( \sum_{i=1}^n [\xitil_i^{(h)} - m_0^{(h)}] \Big)^2}{n^3 \textrm{Log}_2(n)} \Big] \ . \label{eq:vardeclh}
\end{eqnarray}
As to the first summand on the RHS of (\ref{eq:vardecll}), it is enough to notice that $\{(V_0)_{l,l}^{-1}(\xitil_i^{(l)} - m_0^{(l)})^2\}_{i \geq 1}$ is a sequence of i.i.d., $\chi^2(1)$-distributed, real random variables, so that an application of Proposition \ref{prop:teicher} with $r = 4$, in combination with the Lyapunov inequality, yields
$$
\mu_{\theta_0}^{\infty} \Big[ \Big( \sup_{n \geq d+1} \frac{\big{|} \sum_{i=1}^n [(\xitil_i^{(l)} - m_0^{(l)})^2 - (V_0)_{l,l}] \big{|}}{\sqrt{n \textrm{Log}_2(n)}} \Big)^2 \Big] \leq 2(V_0)_{l,l}^2 [\alpha_0(4) + 15^4 \alpha_1(4) ]^{1/2}\ .
$$
For the second summand on the RHS of (\ref{eq:vardecll}), observe that
$$
\frac{\left(\sum_{i=1}^n [\xitil_i^{(l)} - m_0^{(l)}]\right)^4}{n^3 \textrm{Log}_2(n)} = \left(\frac{\big{|} \sum_{i=1}^n [\xitil_i^{(l)} - m_0^{(l)}] \big{|}}{\sqrt{n \textrm{Log}_2(n)}}\right)^4 \frac{[n \textrm{Log}_2(n)]^2}{n^3 \textrm{Log}_2(n)} \leq \left(\frac{\big{|} \sum_{i=1}^n [\xitil_i^{(l)} - m_0^{(l)}] \big{|}}{\sqrt{n \textrm{Log}_2(n)}}\right)^4
$$
yielding, by virtue of Proposition \ref{prop:teicher}, 
$$
\mu_{\theta_0}^{\infty} \Big[ \sup_{n \geq d+1} \frac{\left(\sum_{i=1}^n [\xitil_i^{(l)} - m_0^{(l)}]\right)^4}{n^3 \textrm{Log}_2(n)} \Big] \leq (V_0)_{l,l}^2 [\alpha_0(4) + 3\alpha_1(4)]\ .
$$
To bound the first summand on the RHS of (\ref{eq:vardeclh}), introduce, as a preliminary step, the auxiliary function $\omega_k(\rho) := \ee[(XY - \rho)^k]$, where $k \in \naturals$ and $(X,Y) \in \reals^2$ is a Gaussian random vector with zero means, unit variances and correlation coefficient $\rho$. Then, apply again Proposition \ref{prop:teicher} with $r = 4$, in combination with the Lyapunov inequality, and observe that $1 \leq \omega_2(\rho) = 1 + \rho^2 \leq 2$ 
$\omega_4(\rho) = 9 - 33\rho^2 + 84\rho^2 \leq 60$ lead to
\begin{eqnarray}
&& \mu_{\theta_0}^{\infty} \Big[ \Big( \sup_{n \geq d+1} \frac{\big{|} \sum_{i=1}^n [(\xitil_i^{(l)} - m_0^{(l)})(\xitil_i^{(h)} - m_0^{(h)}) - (V_0)_{l,h}] \big{|}}{\sqrt{n \textrm{Log}_2(n)}} \Big)^2 \Big] \nonumber \\
&\leq& (V_0)_{l,l}(V_0)_{h,h} \omega_2(\rho_{l,h}) \Big[\alpha_0(4) + \alpha_1(4) \Big(\frac{\omega_4(\rho_{l,h})}{\omega_2(\rho_{l,h})} \Big)^4 \Big]^{1/2}
\leq \sqrt{\alpha_0(4) + 60^4\alpha_1(4)} [(V_0)_{l,l}^2 + (V_0)_{h,h}^2] \nonumber
\end{eqnarray}
whre $\rho_{l,h} := (V_0)_{l,h}/\sqrt{(V_0)_{l,l}(V_0)_{h,h}}$. Finally, for the second summand on the RHS of (\ref{eq:vardeclh}), it is enough to notice that
$$
\Big(\sum_{i=1}^n [\xitil_i^{(l)} - m_0^{(l)}]\Big)^2 \cdot \Big(\sum_{i=1}^n [\xitil_i^{(h)} - m_0^{(h)}]\Big)^2 \leq \frac{1}{2} \Big(\sum_{i=1}^n [\xitil_i^{(l)} - m_0^{(l)}]\Big)^4 + \frac{1}{2} \Big(\sum_{i=1}^n [\xitil_i^{(h)} - m_0^{(h)}]\Big)^4\ ,
$$
so that, by resorting to Proposition \ref{prop:teicher} with $r = 4$, one gets
$$
\mu_{\theta_0}^{\infty} \Big[\sup_{n \geq 1} \frac{\left(\sum_{i=1}^n [\xitil_i^{(l)} - m_0^{(l)}]\right)^2 \left(\sum_{i=1}^n [\xitil_i^{(h)} - m_0^{(h)}]\right)^2}{n^3 \textrm{Log}_2(n)} \Big] \leq \frac{1}{2} [(V_0)_{l,l}^2 + (V_0)_{h,h}^2] 
\cdot [\alpha_0(4) + 3\alpha_1(4)]\ .
$$
At this stage, there are all the elements to deduce that
\begin{equation} \label{eq:lastGaussInt}
\mu_{\theta_0}^{\infty}\left[\sup_{n \geq d+1} b_n^2 \text{tr}[V_0 + \hat{V}_n -2(V_0 \hat{V}_n)^{1/2}] \ind\{\text{det}(V_0^{-1}\hat{V}_n) \geq \varepsilon\}\right] \leq c_{\ast} K(\varepsilon, d) \|V_0\|_F
\end{equation}
for a suitable numerical constant $c_{\ast}$ independent of $V_0$ and even of the dimension $d$.


It remains to analyze $\mu_{\theta_0}^{\infty}\big[ \sup_{n \geq d+1} b_n^2 \text{tr}[V_0 + \hat{V}_n -2(V_0 \hat{V}_n)^{1/2}] \ind\{\text{det}(V_0^{-1}\hat{V}_n) < \varepsilon\} \big]$ which, in view of the Boole inequality, can be bounded from above by
$$
\sum_{n \geq d+1} b_n^2 \int_0^{+\infty} \mu_{\theta_0}^{\infty}\left(\left\{ \text{tr}[V_0 + \hat{V}_n -2(V_0 \hat{V}_n)^{1/2}] \ind\{\text{det}(V_0^{-1}\hat{V}_n) < \varepsilon\} > z \right\}\right)\ \ud z\ .
$$
Preliminarily, given any $A > 0$, one gets
$$
\int_0^A \mu_{\theta_0}^{\infty} \left(\left\{ \text{tr}[V_0 + \hat{V}_n -2(V_0 \hat{V}_n)^{1/2}] \ind\{\text{det}(V_0^{-1}\hat{V}_n) < \varepsilon\} > z \right\}\right) \ud z 
\leq A \mu_{\theta_0}^{\infty} \left(\left\{ \text{det}(V_0^{-1}\hat{V}_n) < \varepsilon \right\}\right)
$$
for every $n \in \naturals$. Now, the series $\sum_{n \geq d+1} b_n^2 \mu_{\theta_0}^{\infty} \big( \big\{ \text{det}(V_0^{-1}\hat{V}_n) < \varepsilon \big\} \big)$
turns out to be convergent for any choice of $\varepsilon \in (0,1)$ since, in this case, one can prove that $\mu_{\theta_0}^{\infty} \big( \big\{ \text{det}(V_0^{-1}\hat{V}_n) < \varepsilon \big\} \big)$ goes to zero exponentially fast, with respect to $n$. Indeed, a well-known result concerning the Wishart distribution (see, e.g. Theorem 7.5.3 in \cite{anderson}) states that, for $n > d$, $\text{det}(V_0^{-1}\hat{V}_n)$ has the same distribution (under $\mu_{\theta_0}^{\infty}$) as the product of $d$ independent real random variables, say $Y_1, \dots, Y_d$, where $Y_i$ has a Gamma distribution with scale parameter $\frac{n}{2}$ and shape parameter $\frac{n-i}{2}$. Since $\mu_{\theta_0}^{\infty} 
\big( \big\{ \text{det}(V_0^{-1}\hat{V}_n) < \varepsilon \big\} \big) \leq \sum_{i=1}^d  \mu_{\theta_0}^{\infty} \big( \big\{ Y_i < \varepsilon^{1/d} \big\}\big)$ holds in view of the Boole inequality, then one can resort to the well-known Chernoff bounds to obtain
$$
\mu_{\theta_0}^{\infty} \big( \big\{ Y_i < \varepsilon^{1/d} \big\} \big) \leq \exp\left\{ \frac{n-i}{2} \left[1 - \frac{n \varepsilon^{1/d}}{n-i} + \log\left( \frac{n \varepsilon^{1/d}}{n-i} \right) \right] \right\} 
$$
provided that $\varepsilon^{1/d} < \frac{n-i}{n}$ for any $i \in \{1, \dots, d\}$ and $n \geq d+1$. Since $\sup_{i \leq d, n \geq d+1} \left(\frac{n}{n-i}\right) \leq d+1$, choosing $\varepsilon \leq [2(d+1)]^{-d}$ yields $1 - \frac{n \varepsilon^{1/d}}{n-i} + \log\left( \frac{n \varepsilon^{1/d}}{n-i} \right) \leq \frac12 - \log 2 < 0$  and, hence,
\begin{equation} \label{eq:chernof}
\mu_{\theta_0}^{\infty} \left(\left\{ \text{det}(V_0^{-1}\hat{V}_n) < \varepsilon \right\}\right) \leq d \exp\left\{ \frac{n-d}{2} (1/2 - \log 2) \right\}
\end{equation}
for $n \geq d+1$. Now, one studies $\int_A^{+\infty} \mu_{\theta_0}^{\infty} \Big( \Big\{ \ind\{\text{det}(V_0^{-1}\hat{V}_n) < \varepsilon\} \sum_{j=1}^d [\lambda_j^{(n)} - \log \lambda_j^{(n)} - 1] > z \Big\} \Big) \ud z$
through a splitting of the above integrand into the sum of the following two terms:
$$
\mu_{\theta_0}^{\infty} \Big( \Big\{ \ind\{\text{det}(V_0^{-1}\hat{V}_n) < \varepsilon, \text{tr}(V_0^{-1}\hat{V}_n) > -\eta\log(\text{det}(V_0^{-1}\hat{V}_n))\} \sum_{j=1}^d [\lambda_j^{(n)} - \log \lambda_j^{(n)} - 1] > z/2 \Big\} \Big)
$$
and
$$
\mu_{\theta_0}^{\infty} \Big( \Big\{ \ind\{\text{det}(V_0^{-1}\hat{V}_n) < \varepsilon, \text{tr}(V_0^{-1}\hat{V}_n) \leq -\eta\log(\text{det}(V_0^{-1}\hat{V}_n))\} \sum_{j=1}^d [\lambda_j^{(n)} - \log \lambda_j^{(n)} - 1] > z/2 \Big\} \Big)
$$
with $\eta > 0$. The former can be bounded by $\mu_{\theta_0}^{\infty} \big(\big\{ (1 + 1/\eta) \text{tr}(V_0^{-1}\hat{V}_n) > z/2 \big\}\big)$ while the latter 
can be bounded by $\mu_{\theta_0}^{\infty} \big(\big\{ \text{det}(V_0^{-1}\hat{V}_n) < \varepsilon, \text{det}(V_0^{-1}\hat{V}_n) < \exp\{-\frac{z}{2(1+\eta)}\} \big\}\big)$. At this stage, recall Theorem 7.3.5 in \cite{anderson}, which states that 
$\text{tr}(V_0^{-1}\hat{V}_n))$ has a Gamma distribution with scale parameter $\frac{n}{2}$ and shape parameter $\frac{(n-1)d}{2}$, and apply again the Chernoff bounds to get
$$
\mu_{\theta_0}^{\infty} \Big( \Big\{ \text{tr}(V_0^{-1}\hat{V}_n) > \frac{z\eta}{2(1 + \eta)} \Big\} \Big) \leq \exp\Big\{ \frac{(n-1)d}{2} \Big[1 - \frac{zn\eta}{2(1 + \eta)(n-1)d} + \log\Big(\frac{zn\eta}{2(1 + \eta)(n-1)d}\Big) \Big] \Big\} 
$$
provided that $\frac{z\eta}{2(1 + \eta)} > \frac{(n-1)d}{n}$ holds for any $z \geq A$ and $n \geq d+1$. In particular, this condition is in force if $A \geq 4d(1 + 1/\eta)$. In such a case, one has
$$
\Big[1 - \frac{zn\eta}{2(1 + \eta)(n-1)d} + \log\Big(\frac{zn\eta}{2(1 + \eta)(n-1)d}\Big) \Big] \leq - \frac{(e-2)\eta}{4de(1+\eta)} z
$$
for any $z \geq A$ and $n \geq d+1$. Whence,
\begin{equation} \label{eq:chernoffExt1}
\int_A^{+\infty} \mu_{\theta_0}^{\infty} \Big( \Big\{ \text{tr}(V_0^{-1}\hat{V}_n) > \frac{z\eta}{2(1 + \eta)} \Big\} \Big)\ud z \leq \frac{8e(1+\eta)}{\eta(e-2)(n-1)} \exp\Big\{ - \frac{\eta(e-2) A}{8e(1+\eta)} (n-1) \Big\}\ .
\end{equation}
Apropos of $\mu_{\theta_0}^{\infty} \big(\big\{ \text{det}(V_0^{-1}\hat{V}_n) < \varepsilon, \text{det}(V_0^{-1}\hat{V}_n) < \exp\{-\frac{z}{2(1+\eta)}\} \big\}\big)$, taking $A \geq -2(1+\eta)\log\varepsilon$ implies that
the probability at issue coincides with $\mu_{\theta_0}^{\infty} \big(\big\{ \text{det}(V_0^{-1}\hat{V}_n) < \exp\{-\frac{z}{2(1+\eta)}\} \big\}\big)$, for all $z \geq A$. Therefore, the same reasoning that led to (\ref{eq:chernof}) shows that
\begin{gather}
\mu_{\theta_0}^{\infty} \Big(\Big\{ \text{det}(V_0^{-1}\hat{V}_n) < \exp\{-\frac{z}{2(1+\eta)}\} \Big\} \Big) \nonumber \\
\leq d \exp\Big\{ \frac{n-d}{2} \Big[1 + \log(d+1) - (d+1)\exp\{-\frac{z}{2d(1+\eta)}\} - \frac{z}{2d(1+\eta)} \Big] \Big\}\ . \nonumber
\end{gather}
Now, recalling that $A$ has been chosen in such a way that, $A \geq 4d(1 + 1/\eta)$, one gets
$$
1 + \log(d+1) - (d+1)\exp\{-\frac{z}{2d(1+\eta)}\} - \frac{z}{2d(1+\eta)} \leq \left[\frac{(1 + \log(d+1))\eta}{4d(1+\eta)}  - \frac{1}{2d(1+\eta)} \right] z
$$
which, upon choosing $\eta \leq 1/(1 + \log(d+1))$, leads to
\begin{equation} \label{eq:chernofBIS}
\int_A^{+\infty} \mu_{\theta_0}^{\infty} \Big(\Big\{ \text{det}(V_0^{-1}\hat{V}_n) < \exp\{-\frac{z}{2(1+\eta)}\} \Big\} \Big) \ud z \leq \frac{8d^2(1+\eta)}{(n-d)} \exp\Big\{ - \frac{A}{8d(1+\eta)} (n-d) \Big\}\ .
\end{equation}
In conclusion, for any $\varepsilon \leq [2(d+1)]^{-d}$, it is possible to combine (\ref{eq:chernof})-(\ref{eq:chernofBIS}) and optimize, with respect to the choices of $A$ and $\eta$, to deduce the existence of a positive constant 
$H(\varepsilon, d)$ such that
\begin{equation} \label{eq:lastGauss}
\sum_{n \geq d+1} b_n^2 \int_0^{+\infty} \mu_{\theta_0}^{\infty}\left(\left\{ \text{tr}[V_0 + \hat{V}_n -2(V_0 \hat{V}_n)^{1/2}] \ind\{\text{det}(V_0^{-1}\hat{V}_n) < \varepsilon\} > z \right\}\right)\ \ud z \leq H(\varepsilon, d)
\end{equation}
concluding the proof of (\ref{eq:Main1}), along with the determination of $C_2(\pfrak)$ in (\ref{eq:CpGauss}).

As for the validity of (\ref{eq:Main2}), one considers again (\ref{eq:olkin}) and starts by applying the $d$-dimensional version of the Hartman-Wintner LIL (see, e.g., in Theorem 3.1 of \cite{einmahl}) to write
$$
\pfrak_0^{\infty} \big(\big\{ \limsup_{n \rightarrow +\infty} b_n |\mathbf{m}_0 - \hat{\mathbf{m}}_n| = \sqrt{2} \sigma_{max}(V_0) \big\}\big) = 1
$$ 
where $\sigma_{max}^2(V_0)$ stands for the largest eigenvalue of $V_0$. To deal with the latter summand on the right of (\ref{eq:olkin}), one resorts again to the combination of (\ref{eq:TalagrandGauss}) with (\ref{eq:eigenKullback}) to obtain
\begin{eqnarray}
\sqrt{\text{tr}[V_0 + \hat{V}_n -2(V_0 \hat{V}_n)^{1/2}]} \ind\{\text{det}(V_0^{-1}\hat{V}_n) \geq \varepsilon\} &\leq& \sqrt{2K(\varepsilon,d)\|V_0\|_F} \cdot \|V_0^{-1}\hat{V}_n - \mathrm{Id}_d\|_F \nonumber \\
&\leq& \sqrt{2K(\varepsilon,d)\|V_0^{-1}\|_F} \cdot \|\hat{V}_n - V_0\|_F \nonumber
\end{eqnarray}
for any $\varepsilon \in (0,1)$. Then, observe that $b_n\|\hat{V}_n - V_0\|_F \leq \sum_{l,h = 1}^d \frac{\big{|}\sum_{i=1}^n [\tilde{\eta}_i^{(l,h)} - (V_0)_{l,h}] \big{|}}{\sqrt{n\textrm{Log}_2(n)}}$, by elementary algebra. At this stage, after recalling that 
$$
\sum_{i=1}^n [\tilde{\eta}_i^{(l,h)} - (V_0)_{l,h}] = \sum_{i=1}^n [(\xitil_i^{(l)} - m_0^{(l)})(\xitil_i^{(h)} - m_0^{(h)}) - (V_0)_{l,h}] - \frac{1}{n}\Big(\sum_{i=1}^n [\xitil_i^{(l)} - m_0^{(l)}] \Big) \Big(\sum_{i=1}^n [\xitil_i^{(h)} - m_0^{(h)}] \Big)\ ,
$$ 
apply the classical Hartman-Wintner LIL to conclude that
$$
\pfrak_0^{\infty} \Big(\Big\{ \limsup_{n \rightarrow +\infty} \frac{\big{|}\sum_{i=1}^n [\tilde{\eta}_i^{(l,h)} - (V_0)_{l,h}] \big{|}}{\sqrt{n\textrm{Log}_2(n)}} = \sigma_l\sigma_h\sqrt{2\omega_2(\rho_{l,h})}
\Big\}\Big) = 1
$$
since
$$
\pfrak_0^{\infty} \Big(\Big\{ \limsup_{n \rightarrow +\infty} \frac{\big{|} \big(\sum_{i=1}^n [\xitil_i^{(l)} - m_0^{(l)}] \big) \big(\sum_{i=1}^n [\xitil_i^{(h)} - m_0^{(h)}] \big)\big{|}}{n\sqrt{n\textrm{Log}_2(n)}} = 0
\Big\}\Big) = 1\ .
$$
Noticing that $\sum_{l,h = 1}^d \sigma_l\sigma_h\sqrt{2\omega_2(\rho_{l,h})} \leq d\mathrm{tr}[V_0]$ and that $\sqrt{\|V_0^{-1}\|_F} \cdot \mathrm{tr}[V_0] \leq \sqrt{\mathrm{tr}[V_0]}$, one can conclude the proof after taking into account that
$\pfrak_0^{\infty} \big(\big\{ \limsup_{n \rightarrow +\infty} \ind\{\text{det}(V_0^{-1}\hat{V}_n) \geq \varepsilon\} = 1\big\}\big) = 1$, by strong consistency of $\hat{V}_n$.

\subsection{Proof of Theorem \ref{thm:ExpFam}} \label{sect:proofExpFam}

From the definition of the Kullback-Leibler information, one gets
\begin{eqnarray}
\mathrm{K}(\mu_{\hat{\thetab}_n}\ |\ \mu_{\thetab_0}) &=& \hat{\thetab}_n \cdot [V^{-1}(\hat{\thetab}_n) - V^{-1}(\thetab_0)] - [M(V^{-1}(\hat{\thetab}_n)) - M(V^{-1}(\thetab_0))] \nonumber \\
&=& \int_0^1 (1-s)\ ^t (\hat{\thetab}_n - \thetab_0) \textrm{Hess}[I_{\mu}](\thetab_0 + s(\hat{\thetab}_n - \thetab_0))\ (\hat{\thetab}_n - \thetab_0) \ud s\nonumber
\end{eqnarray}
where the second identity follows from the combination of (\ref{eq:Imu}) with Bernstein's representation of the remainder term in Taylor formula. Then, from (\ref{eq:Talagrand}),
\begin{eqnarray}
\ud_{[\Space]}^{(2)}(\mu_{\thetab_0}, \mu_{\hat{\thetab}_n}) &\leq& \sqrt{C_T(\thetab_0)\int_0^1 (1-s)\ ^t (\hat{\thetab}_n - \thetab_0) \textrm{Hess}[I_{\mu}](\thetab_0 + s(\hat{\thetab}_n - \thetab_0))\ (\hat{\thetab}_n - \thetab_0) \ud s} \nonumber \\
&\leq& \sqrt{C_T(\thetab_0)}\ |\hat{\thetab}_n - \thetab_0|\ \Phi(\hat{\thetab}_n) \nonumber
\end{eqnarray}
and hence, in view of the $k$-dimensional version of the Hartman-Wintner LIL,
\begin{equation}
\limsup_{n \rightarrow +\infty} \sqrt{\frac{n}{\mathrm{Log}_2(n)}} \ud_{[\Space]}^{(2)}(\mu_{\thetab_0}, \mu_{\hat{\thetab}_n}) \leq \sqrt{2 C_T(\thetab_0)} \Phi(\thetab_0) \sigma_{max}(\thetab_0)
\end{equation}
where $\sigma_{max}^2(\thetab_0)$ is the largest eigenvalue of $\mathrm{Hess}[M](V^{-1}(\thetab_0))$. See, e.g., Theorem 3.1 in \cite{einmahl}. This proves (\ref{eq:Main2}) with $b_n = (n/\mathrm{Log}_2(n))^{1/2}$ and $Y(\pfrak_0) = \sqrt{2 C_T(\thetab_0)} \Phi(\thetab_0) \sigma_{max}(\thetab_0)$.

To prove (\ref{eq:Main1}), define $\Theta_i := \{\thetab \in \rk\ |\ |\thetab - \thetab_0| \leq \delta(\thetab_0)\}$ and $\Theta_e := \Theta\setminus\Theta_i$, where $\delta(\thetab_0)$ is chosen so that $\Theta_i$ is a proper subset of $\Theta$, and notice that
\begin{equation}\label{eq:kullbackintext}
\sup_{n \geq 1} b_n |\hat{\thetab}_n - \thetab_0|\ \Phi(\hat{\thetab}_n) \leq \sup_{n \geq 1} b_n |\hat{\thetab}_n - \thetab_0|\ \Phi(\hat{\thetab}_n)\ind\{\hat{\thetab}_n \in \Theta_i\} + \sup_{n \geq 1} b_n |\hat{\thetab}_n - \thetab_0|\ \Phi(\hat{\thetab}_n)\ind\{\hat{\thetab}_n \in \Theta_e\}\ .
\end{equation}
As to the first summand on the right-hand side, for any $r > 2$, one can write
$$
\mu_{\thetab_0}^{\infty}\left[\sup_{n \geq 1} b_n |\hat{\thetab}_n - \thetab_0|\ \Phi(\hat{\thetab}_n)\ind\{\hat{\thetab}_n \in \Theta_i\}\right] \leq \max_{\etab \in \Theta_i} \Phi(\etab) \sum_{j=1}^k \left(\mu_{\thetab_0}^{\infty}[\sup_{n \geq 1} b_n^r |\hat{\theta}_n^{(j)} - \theta_0^{(j)}|^r]\right)^{1/r}
$$
where $\hat{\theta}_n^{(j)}$ and $\theta_0^{(j)}$ stand for the $j$-th component of $\hat{\thetab}_n$ and $\thetab_0$, respectively. At this stage, invoke Proposition \ref{prop:teicher} to conclude that
\begin{gather}
\mu_{\thetab_0}^{\infty}\left[\sup_{n \geq 1} b_n |\hat{\thetab}_n - \thetab_0|\ \Phi(\hat{\thetab}_n)\ind\{\hat{\thetab}_n \in \Theta_i\}\right] \\
\leq \max_{\etab \in \Theta_i} \Phi(\etab) \sum_{j=1}^k \sigma_j(\thetab_0) \left[\alpha_0(r) + \alpha_1(r) \left(\frac{\mu_{\thetab_0}^{\infty}[|t^{(j)}(\tilde{\xi}_1) - \theta_0^{(j)}|^r]}{\sigma_j(\thetab_0)^r}\right)^{\lceil r \rceil}\right]^{1/r}
\end{gather}
where $\sigma_j(\thetab_0) := \sqrt{\partial_{j,j}^2 M(V^{-1}(\thetab_0))}$ and $t^{(j)}(\tilde{\xi}_1)$ stands for the $j$-th component of $\tb(\tilde{\xi}_1)$. This upper bound represents a first contribution to the determination of $C(\pfrak_0)$, to be now completed by bounding the second summand on the right-hand side of (\ref{eq:kullbackintext}). Apropos of this, one can start with the following general considerations, valid for any $n_0 \in \naturals$:
\begin{eqnarray}
&& \mu_{\thetab_0}^{\infty}\left[\sup_{n \geq n_0} b_n |\hat{\thetab}_n - \thetab_0|\ \Phi(\hat{\thetab}_n)\ind\{\hat{\thetab}_n \in \Theta_e\}\right]
\nonumber \\
&=& \int_0^{+\infty} \mu_{\thetab_0}^{\infty}\left[\bigcup_{n \geq n_0} \left\{ b_n |\hat{\thetab}_n - \thetab_0|\ \Phi(\hat{\thetab}_n)\ind\{\hat{\thetab}_n \in \Theta_e\} > t\right\}\right] \ud t \nonumber \\
&\leq& \sum_{n \geq n_0} b_n \int_0^{\tau(\thetab_0)} \mu_{\thetab_0}^{\infty}\left[|\hat{\thetab}_n - \thetab_0|\ \Phi(\hat{\thetab}_n)\ind\{\hat{\thetab}_n \in \Theta_e\} > t \right] \ud t \label{eq:boole1} \\
&+& \sum_{n \geq n_0} b_n \int_{\tau(\thetab_0)}^{+\infty} \mu_{\thetab_0}^{\infty}\left[|\hat{\thetab}_n - \thetab_0|\ \Phi(\hat{\thetab}_n)\ind\{\hat{\thetab}_n \in \Theta_e\} > t \right] \ud t\ . \label{eq:boole2}
\end{eqnarray}
The term in (\ref{eq:boole1}) can be bounded thanks to the fact that $\left\{|\hat{\thetab}_n - \thetab_0|\ \Phi(\hat{\thetab}_n)\ind\{\hat{\thetab}_n \in \Theta_e\} > t\right\} \subseteq \{\hat{\thetab}_n \in \Theta_e\}$, yielding
$$
\sum_{n \geq n_0} b_n \int_0^{\tau(\thetab_0)} \mu_{\thetab_0}^{\infty}\left[|\hat{\thetab}_n - \thetab_0|\ \Phi(\hat{\thetab}_n)\ind\{\hat{\thetab}_n \in \Theta_e\} > t \right] \ud t \leq \tau(\thetab_0) \sum_{n \geq n_0} b_n \mu_{\thetab_0}^{\infty}\left[|\hat{\thetab}_n - \thetab_0| > \delta(\thetab_0)\right]\ .
$$
Here and in other points of the present proof, some results in \cite{yurinskii} are invoked to estimate tail probabilities of the type of $\mu_{\thetab_0}^{\infty}\left[|\hat{\thetab}_n - \thetab_0| > a\right]$. More precisely, the analysis of (\ref{eq:boole1}) is based on the Corollary on page 491 of that paper, whose applicability to the present context follows after checking (2.1) therein. Since vectors $\tb(\tilde{\xi}_j) - \thetab_0$ correspond to vectors $\boldsymbol{\xi}_j$ in \cite{yurinskii}, in the place of $\ee[|\boldsymbol{\xi}_j|^m]$ therein one considers
$$
\mu_{\thetab_0}^{\infty}[|\tb(\tilde{\xi}_j) - \thetab_0|^m] \leq k^{m/2 - 1} \sum_{i=1}^k \mu_{\thetab_0}^{\infty}[|t^{(i)}(\tilde{\xi}_j) - \theta_0^{(i)}|^m]\ .
$$
Thanks to the existence of the moment generating function of $\tb(\tilde{\xi}_j) - \thetab_0$, the classical Cauchy estimate shows that
$$
\mu_{\thetab_0}^{\infty}[|t^{(i)}(\tilde{\xi}_j) - \theta_0^{(i)}|^m] \leq m! \frac{C(r(\thetab_0))}{r(\thetab_0)^m}
$$
is valid for all $m \in \naturals$ with a suitable chosen $r(\thetab_0) > 0$, provided that $C(r(\thetab_0))$ is the maximum (as $i = 1, \dots, k$) of the maximum modulus of the moment generating function of $|t^{(i)}(\tilde{\xi}_j) - \theta_0^{(i)}|$ on $\{z \in \mathbb{C}\ |\ |z| = r(\thetab_0)\}$. Plainly, the last inequality entails (2.1) in \cite{yurinskii} with specific $H = H(\thetab_0)$ and $b_j^2 = \mu_{\thetab_0}^{\infty}[|\tb(\tilde{\xi}_j) - \thetab_0|^2] =: B(\thetab_0)^2$. Putting $B_n^2 := nB^2(\thetab_0)$ and $x := \sqrt{n}\delta(\thetab_0)/B(\thetab_0)$, the thesis of the aforesaid corollary yields
\begin{eqnarray}
\sum_{n \geq n_0} b_n \mu_{\thetab_0}^{\infty}\left[|\hat{\thetab}_n - \thetab_0| > \delta(\thetab_0)\right] &\leq& \sum_{n \geq n_0} b_n \exp\left\{
\frac{-n \delta(\thetab_0)^2}{B(\thetab_0)^2 + 1,62\delta(\thetab_0)H(\thetab_0)}\right\} \nonumber \\
&\leq& c \left(\frac{\delta(\thetab_0)^2}{B(\thetab_0)^2 + 1,62\delta(\thetab_0)H(\thetab_0)}\right)^{-2}
\exp\left\{\frac{-n_0\delta(\thetab_0)^2}{B(\thetab_0)^2 + 1,62\delta(\thetab_0)H(\thetab_0)}\right\} \nonumber \\
\label{eq:yurinski1}
\end{eqnarray}
as a consequence of the elementary inequality
\begin{equation} \label{eq:serie}
\sum_{n \geq n_0} b_n e^{-nz} \leq \frac{c}{z^2}e^{-n_0 z}
\end{equation}
valid for $z > 0$ with some constant $c>0$. Estimate (\ref{eq:yurinski1}) provides the required upper bound for (\ref{eq:boole1}). As for (\ref{eq:boole2}), start by writing
\begin{eqnarray}
&& \mu_{\thetab_0}^{\infty}[|\hat{\thetab}_n - \thetab_0| \Phi(\hat{\thetab}_n) > t] \nonumber \\
&\leq& \mu_{\thetab_0}^{\infty}[|\hat{\thetab}_n - \thetab_0| \Phi(\hat{\thetab}_n) > t, \Phi(\hat{\thetab}_n) \leq \sigma(t)] + \mu_{\thetab_0}^{\infty}[|\hat{\thetab}_n - \thetab_0| \Phi(\hat{\thetab}_n) > t, \Phi(\hat{\thetab}_n) > \sigma(t)] \nonumber \\
&\leq& \mu_{\thetab_0}^{\infty}[|\hat{\thetab}_n - \thetab_0| > t/\sigma(t)] + \mu_{\thetab_0}^{\infty}[\Phi(\hat{\thetab}_n) > \sigma(t)] \nonumber \\
&\leq& \mu_{\thetab_0}^{\infty}[|\hat{\thetab}_n - \thetab_0| > t/\sigma(t)] + \mu_{\thetab_0}^{\infty}[\Phi(\hat{\thetab}_n) > \sigma(t), |\hat{\thetab}_n - \thetab_0| \leq \rho(t)] + \mu_{\thetab_0}^{\infty}[|\hat{\thetab}_n - \thetab_0| > \rho(t)] \nonumber \\
&\leq& \mu_{\thetab_0}^{\infty}[|\hat{\thetab}_n - \thetab_0| > t/\sigma(t)] + \mu_{\thetab_0}^{\infty}[\Phi_1(|\hat{\thetab}_n - \thetab_0|) > \frac{\sigma(t)}{2}] + \mu_{\thetab_0}^{\infty}[\Phi_2(\hat{\thetab}_n) > \frac{\sigma(t)}{2}, |\hat{\thetab}_n - \thetab_0| \leq \rho(t)] \nonumber \\
&+& \mu_{\thetab_0}^{\infty}[|\hat{\thetab}_n - \thetab_0| > \rho(t)] \nonumber
\end{eqnarray}
and, after introducing the function $m(t) := \min\{\rho(t), \frac{t}{\sigma(t)}, \Phi_1^{-1}(\frac{\sigma(t)}{2})\}$, the sum of the first, second and fourth term can be bounded by
$$
3\exp\left\{\frac{-n m(t)^2}{B(\thetab_0)^2 + 1,62m(t)H(\thetab_0)}\right\} \leq 3\exp\left\{\frac{-n m(t)}{2H(\thetab_0)}\right\}\ \ \ \ \ \ (t > \tau(\thetab_0))
$$
invoking once again the Corollary in \cite{yurinskii} with $x := \sqrt{n}m(t)/B(\thetab_0)$. The validity of the above inequality rests on the fact that, by properly choosing $\rho$ and $\sigma$, $m$ turns out to be monotone and diverging at infinity. Now, from (\ref{eq:serie}),
$$
\sum_{n \geq n_0} 3b_n\exp\left\{\frac{-n m(t)}{2H(\thetab_0)}\right\} \leq 3c \left(\frac{m(t)}{2H(\thetab_0)}\right)^{-2} \exp\left\{\frac{- n_0 m(t)}{2H(\thetab_0)}\right\}\ \ \ \ \ \ (t > \tau(\thetab_0))
$$
and one finds that the right-hand side is integrable on $[\tau(\thetab_0), +\infty)$, thanks to (\ref{eq:CostExp}), after a proper choice of $n_0$. At this stage, it remains only to bound $\mu_{\thetab_0}^{\infty}[\Phi_2(\hat{\thetab}_n) > \frac{\sigma(t)}{2}, |\hat{\thetab}_n - \thetab_0| \leq \rho(t)] = \mu_{\thetab_0}^{\infty}[\hat{\thetab}_n \in R_1(t)]$, a task which will be carried out by an application of some large deviation estimate. In fact, recall that
$$
\mu_{\thetab_0}^{\infty}[\hat{\thetab}_n \in R] \leq e^{-n \inf_{\thetab \in R} I_{\thetab_0}(\thetab)}
$$
holds true for every closed convex subset $R$ of $\Theta$. See, e.g., Section 2 of \cite{borovkov}. Therefore, since $R_1(t)$ is compact, consider a covering $\{V_i(t)\}_{i=1, \dots,N(t)}$ of $R_1(t)$, made by closed and convex subsets of $\Theta$, which is entirely contained in $R_2$, to obtain
$$
\mu_{\thetab_0}^{\infty}[\hat{\thetab}_n \in R_1(t)] \leq \mu_{\thetab_0}^{\infty}[\hat{\thetab}_n \in \cup_{i=1}^{N(t)} V_i(t)] \leq \sum_{i=1}^{N(t)} \mu_{\thetab_0}^{\infty}[\hat{\thetab}_n \in V_i(t)] \leq N(t) e^{-n h(t)}
$$
where $h(t) := \inf\{I_{\thetab_0}(\thetab)\ |\ \Phi_2(\thetab) \geq \frac{\sigma(t)}{4}, |\thetab - \thetab_0| \leq \rho(t)\}$ and $N(t)$ is an abbreviation for $N_{\Phi_2}(\rho(t), \sigma(t))$. The conclusion is reached by first applying (\ref{eq:serie}), giving
$$
\sum_{n \geq n_0} b_n e^{-n h(t)} \leq c\left(\frac{1}{h(t)}\right)^2 e^{-n_0 h(t)}
$$
and then invoking (\ref{eq:CostExp}) to prove the integrability of $N(t) \left(\frac{1}{h(t)}\right)^2 e^{-n_0 h(t)}$ on $[\tau(\thetab_0), +\infty)$. After choosing $n_0$ properly, sums of the type $\sum_{n =1}^{n_0-1}$ admits trivial bounds which, combined with the rest of the proof, lead to the desired result.

\subsection{Proof of Theorem \ref{thm:BD}} \label{sect:proofBD}

First, observe that the following two identities
$$
\left[\ud_{[\pms_p]}^{(p)}\left(\pi(\xisun), \delta_{\hat{\pfrak}_n}\right)\right]^p = \int_{\pms} \left(\ud_{\pms}^{(p)}(\pfrak, \hat{\pfrak}_n)\right)^p \pi(\xisun, \ud \pfrak) = \rho\left[\left(\ud_{\pms}^{(p)}(\randommeasure, \hat{\pfrak}_n)\right)^p \ |\ \xinsu\right]
$$
hold with $\rho$-probability 1, the former being a direct consequence of the definition (\ref{eq:Wp}). To derive the latter, suffice it to apply the disintegration theorem (see, e.g., Theorem 6.4 in \cite{ka}) after recalling that 
$\hat{\pfrak}_n : (\Omega, \mathscr{F}_n) \rightarrow (\Space, \mathscr{B}(\pms_p))$ is measurable. In the next step, parallelling the reasoning in \cite{blackdub} to prove Theorem 2 therein, one gets
\begin{eqnarray}
\limsup_{n \rightarrow \infty} \rho\left[\left(b_n \ud_{\pms}^{(p)}(\randommeasure, \hat{\pfrak}_n)\right)^p \ |\ \xinsu\right]
&\leq& \lim_{k \rightarrow +\infty} \sup_{n \geq k, h \geq k} \rho\left[\left(b_h \ud_{\pms}^{(p)}(\randommeasure, \hat{\pfrak}_h)\right)^p \ |\ \xinsu\right] \nonumber \\
&\leq& \lim_{k \rightarrow +\infty} \sup_{n \geq k} \rho\left[\sup_{h \geq m} \left(b_h \ud_{\pms}^{(p)}(\randommeasure, \hat{\pfrak}_h)\right)^p \ |\ \xinsu\right] \nonumber
\end{eqnarray}
for every $m \in \naturals$. Then, one has that $X_m := \sup_{h \geq m} \left(b_h \ud_{\pms}^{(p)}(\randommeasure, \hat{\pfrak}_h)\right)^p$ belongs to $\mathrm{L}^1(\Omega, \mathscr{F}, \rho)$, since
$$
\rho[X_m] \leq \rho[X_1] = \rho\left[\rho[X_1\ |\ \randommeasure]\right] \leq \rho[C_p(\randommeasure)] < +\infty\ .
$$
Hence, the martingale convergence theorem gives
$$
\lim_{k \rightarrow +\infty} \sup_{n \geq k} \rho\big[X_m \ |\ \xinsu\big] = \rho \Big[ X_m \ \Big{|}\ \sigma\Big{(}\bigcup_{n \geq 1}\sigma(\xinsu)\Big{)}\Big]
$$
with $\rho$-probability 1, and now it is enough to observe that $X_m$ is a $\sigma\left(\bigcup_{n \geq 1} \mathscr{F}_n\right)$-measurable real random variable, for any $m \in \naturals$, since $\rho[\empiric \Rightarrow \randommeasure, \ \text{as}\ n \rightarrow +\infty] = 1$.
Whence,
$$
\rho\left[\left\{\lim_{k \rightarrow +\infty} \sup_{n \geq k} \rho\left[X_m \ |\ \xinsu\right] = X_m \right\}\right] = 1
$$
for every $m \in \naturals$, and $\limsup_{n \rightarrow \infty} \left[\ud_{[\pms_p]}^{(p)}\left(\pi(\xisun), \delta_{\hat{\pfrak}_n}\right)\right]^p \leq \limsup_{m \rightarrow +\infty} X_m \leq [Y_p(\randommeasure)]^p$ hold
with $\rho$-probability 1, yielding (\ref{eq:MainP}).

\subsection{Proof of Theorem \ref{thm:BDpred}} \label{sect:proofBDpred}

First, observe that (\ref{eq:MainPred}) follows immediately from the combination of (\ref{eq:savare1}) with (\ref{eq:MainP}). Thus, given $p \geq 1$, one proves (\ref{eq:savare1}) as a consequence of the following more general statement: for $\pfrak_0 \in \pms_p$ and $\zeta \in [\pms_p]$ there holds
\begin{equation} \label{eq:BirkhoffInfinity}
\ud^{(p)}_{[\pms_p]}\left(\rho[\zeta] \circ \empiricm^{-1}, \pfrak_0^{\infty} \circ \empiricm^{-1} \right) \leq \ud^{(p)}_{[\pms_p]}\left(\zeta, \delta_{\pfrak_0}\right)
\end{equation}
for every $m \in \naturals$, $\rho[\zeta]$ being the p.m. on $(\Omega, \mathscr{F}) = (\Space^{\infty}, \mathscr{B}(\Space^{\infty}))$ given by $\rho[\zeta](\cdot) := \int_{\pms_p} \pfrak^{\infty}(\cdot) \zeta(\ud\pfrak)$. Now, for simplicity, assume that $\zeta = \sum_{j=1}^N \alpha_j \delta_{\pfrak_j}$, where $\pfrak_1, \dots, \pfrak_N \in \pms_p$ and $\alpha_1, \dots, \alpha_N \in [0,1]$ with $\sum_{j=1}^N \alpha_j = 1$, and define $\overline{\eta}[\pfrak_0, \pfrak_j] \in \mathcal{F}(\pfrak_0, \pfrak_j)$ to be an optimal coupling for $\ud^{(p)}_{\pms}$, i.e.\! $\ud^{(p)}_{\pms}(\pfrak_0, \pfrak_j) = \big{(} \int_{\Space^2} \ud_{\Space}(x,y)^p\ \overline{\eta}[\pfrak_0, \pfrak_j](\ud x \ud y) \big{)}^{1/p}$, for $j = 1, \dots, N$. Then, put $\overline{\gamma}_m(A_1 \times \dots \times A_m \times B_1\times \dots \times B_m ) := \sum_{j=1}^N \alpha_j \{\prod_{i=1}^m \overline{\eta}[\pfrak_0, \pfrak_j](A_i \times B_i)\}$ for $A_1, \dots, A_m, B_1,\dots, B_m \in \mathscr{B}(\Space)$, noticing that $\overline{\gamma}_m \in \mathcal{F}(\pfrak_0^m, \rho_m[\zeta])$, where $\rho_m[\zeta](\cdot) := \int_{\pms_p} \pfrak^m(\cdot) \zeta(\ud\pfrak)$. Whence, 
from the definition of $\ud^{(p)}_{[\pms_p]}$,
$$
\ud^{(p)}_{[\pms_p]}\Big(\rho[\zeta] \circ \empiricm^{-1}, \pfrak_0^{\infty} \circ \empiricm^{-1} \Big) \leq \Big[ \int_{\Space^{2m}} \ud^{(p)}_{\pms}\Big(\frac{1}{m} \sum_{i=1}^m \delta_{x_i}, \frac{1}{m} \sum_{i=1}^m \delta_{y_i}\Big)^p \overline{\gamma}_m(\ud x_1 \dots \ud x_m \ud y_1 \dots \ud y_m) \Big]^{1/p}\ .
$$
Since a combination of the convexity of $\ud^{(p)}_{\pms}$ with the Lyapunov inequality for moments yields $\ud^{(p)}_{\pms}\left(\frac{1}{m} \sum_{i=1}^m \delta_{x_i}, \frac{1}{m} \sum_{i=1}^m \delta_{y_i}\right)^p \leq \frac{1}{m} 
\sum_{i=1}^m \ud^{(p)}_{\pms}(\delta_{x_i}, \delta_{y_i})^p = \frac{1}{m} \sum_{i=1}^m \ud_{\Space}(x_i, y_i)^p$, one gets
\begin{eqnarray}
\ud^{(p)}_{[\pms_p]}\left(\rho[\zeta] \circ \empiricm^{-1}, \pfrak_0^{\infty} \circ \empiricm^{-1} \right) &\leq& \Big[ \int_{\Space^{2m}} \Big\{ \frac{1}{m} \sum_{i=1}^m \ud_{\Space}(x_i, y_i)^p \Big\} 
\overline{\gamma}_m(\ud x_1 \dots \ud x_m \ud y_1 \dots \ud y_m) \Big]^{1/p} \nonumber \\
&=& \Big[ \sum_{j=1}^N \alpha_j \int_{\Space^2} \ud_{\Space}(x, y)^p \overline{\eta}[\pfrak_0, \pfrak_j](\ud x \ud y) \Big]^{1/p} \!\! \!\! = \Big[ \sum_{j=1}^N \alpha_j \ud^{(p)}_{\pms}(\pfrak_0, \pfrak_j)^p \Big]^{1/p}\!\! . \nonumber
\end{eqnarray}
Hence, the identities $\sum_{j=1}^N \alpha_j \ud^{(p)}_{\pms}(\pfrak_0, \pfrak_j)^p = \int_{\pms_p} \ud^{(p)}_{\pms}(\pfrak_0, \pfrak)^p \zeta(\ud \pfrak) = \ud^{(p)}_{[\pms_p]}\left(\zeta, \delta_{\pfrak_0}\right)^p$ lead to the proof of 
(\ref{eq:BirkhoffInfinity}), under the assumption that $\zeta = \sum_{j=1}^N \alpha_j \delta_{\pfrak_j}$.

To complete the argument, invoke (the proof of) Lemma 11.8.4 in \cite{du} to show that any p.m.\! $\zeta \in [\pms_p]$ admits an approximating sequence $\{\zeta^{(N)}\}_{N \geq 1}$ with $\zeta^{(N)} = \sum_{j=1}^N \alpha_j^{(N)} \delta_{\pfrak_j^{(N)}}$, such that $\ud^{(p)}_{[\pms_p]}\left(\zeta, \zeta^{(N)}\right) \rightarrow 0$ as $N \rightarrow \infty$. Now, $\rho[\zeta^{(N)}] \circ \empiricm^{-1} \Rightarrow \rho[\zeta] \circ \empiricm^{-1}$ holds for every $m \in \naturals$ by virtue of the continuous mapping theorem (cf., e.g., Theorems 4.27 and 4.29 in \cite{ka}). Hence, after noticing that the sequence 
$\big\{\int_{\pms_p} \int_{\Space} \ud_{\Space}(x_0, x)^p \pfrak(\ud x) \zeta^{(N)}(\ud\pfrak)\big\}_{N \geq 1}$ is bounded, in view of Lemma 11.8.4 in \cite{du},
 exploit the lower semicontinuity of $\ud^{(p)}_{[\pms_p]}$ (see Proposition 7.1.3 in \cite{amgisa}), i.e.\! $\ud^{(p)}_{[\pms_p]}\left(\rho[\zeta] \circ \empiricm^{-1}, \pfrak_0^{\infty} \circ \empiricm^{-1} \right) \leq \liminf_{N \rightarrow \infty} \ud^{(p)}_{[\pms_p]}\left(\rho[\zeta^{(N)}] \circ \empiricm^{-1}, \pfrak_0^{\infty} \circ \empiricm^{-1} \right)$, to deduce the validity of (\ref{eq:BirkhoffInfinity}) for any $\zeta \in [\pms_p]$.

\vspace{1cm}



\begin{multicols}{2}
\footnotesize{\textsc{emanuele dolera} (corresponding author) \\
\textsc{dipartimento di matematica ``felice casorati'' \\
universit\`a degli studi di pavia \\
via ferrata 1, 27100 pavia, italy \\
e-mail:} emanuele.dolera@unipv.it}

\columnbreak

\footnotesize{\textsc{eugenio regazzini \\
dipartimento di matematica ``felice casorati'' \\
universit\`a degli studi di pavia \\
via ferrata 1, 27100 pavia, italy \\
e-mail:} eugenio.regazzini@unipv.it}

\end{multicols}


\begin{thebibliography}{99}

\bibitem{ajtai} \textsc{Ajtai, M., Koml\'{o}s, J.} and  \textsc{Tusn\'{a}dy, G.} (1984). On optimal matchings. \emph{Combinatorica} $\mathbf{4}$, 259-264.

\bibitem{ald} \textsc{Aldous, D.J.} (1985). Exchangeability and related topics. \emph{Ecole d'Et\'{e} de Probabilit\'{e}s de Saint-Flour, XIII} 1-198. Lecture Notes in Math. $\mathbf{1117}$ Springer, Berlin.

\bibitem{amgisa} \textsc{Ambrosio, L.}, \textsc{Gigli, N.} and \textsc{Savar\'{e}, G.} (2008). \emph{Gradient Flows in Metric Spaces and in the Space of Probability Measures}. $2^{nd}$ ed. Birkh\"{a}user, Basel.

\bibitem{ambrosio} \textsc{Ambrosio, L., Stra, F.} and \textsc{Trevisan, D.} (2016). A PDE approach to a 2-dimensional matching problem. \emph{arXiv}:1611.04960.

\bibitem{anderson} \textsc{Anderson, T.W.} (1984). \emph{An Introduction to Multivariate Statistical Analysis}. $2^{nd}$ ed. Wiley, New York.

\bibitem{barndorff} \textsc{Barndorff-Nielsen, 0.} (1978). \emph{Information and Exponential Families in Statistical Theory}. Wiley, New York.

\bibitem{BCPR} \textsc{Berti, P., Crimaldi, I., Pratelli, L.} and \textsc{Rigo, P.} (2009). Rate of convergence of predictive distributions for dependent data. 
\emph{Bernoulli} $\mathbf{15}$, 1351-1367.

\bibitem{BPR} \textsc{Berti, P., Pratelli, L.} and \textsc{Rigo, P.} (2012). Limit theorems for empirical processes based on dependent data.
\emph{Electron. J. Probab.} $\mathbf{17}$, 1-18.

\bibitem{blackdub} \textsc{Blackwell, D.} and \textsc{Dubins, L.E.} (1962). Merging of opinions with increasing information.
\emph{Ann. Math. Statist.} $\mathbf{33}$, 882-886.

\bibitem{BobkovLedoux} \textsc{Bobkov, S.} and \textsc{Ledoux, M.} (2016). One-dimensional empirical measures, order statistics, and Kantorovich transport distances. To appear in: \emph{Memoirs of the AMS}.

\bibitem{boissard} \textsc{Boissard, P.} (2011). Simple bounds for the convergence of empirical and occupation measures in 1-Wasserstein distance. \emph{Electron. J. Probab.} $\mathbf{16}$, 2296-2333.

\bibitem{boissardlegouic} \textsc{Boissard, P.} and \textsc{Le Gouic, T.} (2014). On the mean speed of convergence of empirical and occupation
measures in Wasserstein distance. \emph{Ann. Inst. Henri Poincar\'e--Probab. Stat.} $\mathbf{50}$, 539-563.

\bibitem{bolley} \textsc{Bolley, F., Guillin, A.} and \textsc{Villani, C.} (2007). Quantitative concentration inequalities for empirical measures on non-compact spaces. \emph{Probab. Theory Related Fields} $\mathbf{137}$, 541-593.

\bibitem{borovkov} \textsc{Borovkov, A.A.} and \textsc{Mogulskii, A.A.} (2012). Chebyshev-type exponential inequalities. \emph{Theory Probab. Appl.} $\mathbf{56}$, 21-43.

\bibitem{brown} \textsc{Brown, L.D.} (1986). \emph{Fundamentals of Statistical Exponential Families with Application in Statistical Decision Theory}. Institute of Mathematical Statistics Lecture Notes-Monograph Series, 9. Hayward, California.

\bibitem{cantelli} \textsc{Cantelli, F.P.} (1933). Sulla determinazione empirica delle leggi di probabilit\`a.
\emph{Giorn. Ist. Ital. Attuari} $\mathbf{4}$, 421-424.

\bibitem{caracciolo} \textsc{Caracciolo, S., Lucibello, G., Parisi, G.,} and \textsc{Sicuro, G.} (2014). Scaling hypothesis for the Euclidean bipartite matching problem. \emph{Phys. Rev. E} $\mathbf{90}$, 012118.

\bibitem{chte} \textsc{Chow, Y. S.} and \textsc{Teicher, H} (1997). \emph{Probability Theory. Independence, Interchangeability, Martingales}. $3^{rd}$ ed. Springer, New York.

\bibitem{cifdorega} \textsc{Cifarelli, D.M., Dolera, E.} and \textsc{Regazzini, E.} (2016). Frequentistic approximations to Bayesian prevision of exchangeable random elements. \emph{Internat. J. Approx. Reason.} $\mathbf{78}$, 138-152.

\bibitem{cifdoregaNOTE} \textsc{Cifarelli, D.M., Dolera, E.} and \textsc{Regazzini, E.} (2017). Note on ``Frequentistic approximation to Bayesian prevision of exchangeable random elements''. \emph{Internat. J. Approx. Reason.} $\mathbf{86}$, 26-27.

\bibitem{cox} \textsc{Cox, D.R.} (1975). Partial likelihood. \emph{Biometrika} $\textbf{62}$, 269-276.

\bibitem{dereich} \textsc{Dereich, S., Scheutzow, M.} and \textsc{Schottstedt, R.} (2013). Constructive quantization: approximation by empirical measures. \emph{Ann. Inst. Henri Poincar\'{e} Probab. Stat.} $\mathbf{49}$, 1183-1203.

\bibitem{lugosi} \textsc{Devroye, L., Gy\"orfi, L.} and \textsc{Lugosi, G.} (1996). \emph{A Probabilistic Theory of Pattern Recognition}. Springer, New York.

\bibitem{diafree1} \textsc{Diaconis, P.} and \textsc{Freedman, D.} (1986). On consistency of Bayes estimates. \emph{Ann. Statist.} $\mathbf{14}$, 1-26.

\bibitem{dowson} \textsc{Dowson, D.C.} and \textsc{Landau, B.V.} (1982). The Fr\'{e}chet distance between multivariate normal distributions. \emph{J. Multivariate Anal.} $\mathbf{12}$, 450-455.

\bibitem{dudleyGC} \textsc{Dudley, R.M.} (1969). The speed of mean Glivenko-Cantelli convergence. \emph{Ann. Math. Statist.} $\mathbf{40}$, 40-50.

\bibitem{dudleyCLT} \textsc{Dudley, R.M.} (1999). \emph{Uniform Central Limit Theorems}. Cambridge University Press, Cambridge.

\bibitem{du} \textsc{Dudley, R.M.} (2002). \emph{Real Analysis and Probability}. Cambridge University Press, Cambridge.

\bibitem{efron} \textsc{Efron, B.} (1979). Bootstrap methods: another look at the jackknife. \emph{Ann. Statist.} $\mathbf{7}$, 1-26.

\bibitem{efronBOOK} \textsc{Efron, B.} (2010). \emph{Large-Scale Inference}. Cambridge University Press, Cambridge.

\bibitem{einmahl} \textsc{Einmahl, U.} (2016). Law of the iterated logarithm type results for random vectors with infinite second moments. \emph{Math. Appl.}  $\mathbf{44}$, 167-181.

\bibitem{definetti30} de \textsc{Finetti, B.} (1930). Funzione caratteristica di un fenomeno aleatorio. \emph{Atti Reale Accademia Nazionale dei Lincei, Mem.} $\textbf{4}$, 86-133.

\bibitem{definetti33} de \textsc{Finetti, B.} (1933). Sull'approssimazione empirica di una legge di probabilit\`a. \emph{Giorn. Istit. Ital. Attuari} $\textbf{4}$, 415-420.

\bibitem{definettiLincei} de \textsc{Finetti, B.} (1933). La legge dei grandi numeri nel caso dei numeri aleatori equivalenti. \emph{Atti Reale Accademia Nazionale dei Lincei, Rend.} $\textbf{18}$, 203-207.

\bibitem{definetti37} de \textsc{Finetti, B.} (1937). La pr\'{e}vision: ses lois logiques, ses sources subjectives. \emph{Ann. Inst. H. Poincar\'{e}} $\textbf{7}$, 1-68.

\bibitem{fourn} \textsc{Fournier, N.} and \textsc{Guillin, A.} (2015). On the rate of convergence in Wasserstein distance of the empirical measure. \emph{Probab. Theory Related Fields} $\mathbf{162}$, 707-738.

\bibitem{glivenko} \textsc{Glivenko, V.I.} (1933). Sulla determinazione empirica delle leggi di probabilit\`a.
\emph{Giorn. Istit. Ital. Attuari} $\mathbf{4}$, 92-99.

\bibitem{gozlanAOP} \textsc{Gozlan, N.} (2009). A characterization of dimension free concentration in terms of transport inequalities. \emph{Ann. Probab.} $\mathbf{37}$, 2480-2498.

\bibitem{gozlan1} \textsc{Gozlan, N.} and \textsc{L\'{e}onard, C.} (2007). A large deviation approach to some transportation cost inequalities. \emph{Probab. Theory Related Fields} $\mathbf{139}$, 235-283.

\bibitem{gozlanSurvey} \textsc{Gozlan, N.} and \textsc{L\'{e}onard, C.} (2010). Transport inequalities. A survey. \emph{Markov Processes Related Fields} $\mathbf{16}$, 635-736.

\bibitem{horn} \textsc{Horn, R.A.} and \textsc{Johnson, C.R.} (2013). \emph{Matrix Analysis}. $2^{nd}$ ed. Cambridge University Press, Cambridge.

\bibitem{horowitz} \textsc{Horowitz, J.} and \textsc{Karandikar, R.L.} (1994). Mean rates of convergence of empirical measures in the Wasserstein metric. \emph{J. Comput. Appl. Math.} $\mathbf{55}$, 261-273.

\bibitem{ka} \textsc{Kallenberg, O.} (2002). \emph{Foundations of Modern Probability}. $2^{nd}$ ed. Springer-Verlag, New York.

\bibitem{kolmoGC} \textsc{Kolmogorov, A.N.} (1933). Sulla determinazione empirica di una legge di distribuzione.
\emph{Giorn. Istit. Ital. Attuari} $\mathbf{4}$, 83-91.

\bibitem{maritz} \textsc{Maritz, J.S.} and \textsc{Lwin, T.} (1989). \emph{Empirical Bayes Methods}. $2^{nd}$ ed. Chapman and Hall, London.

\bibitem{massart} \textsc{Massart, J.E.} (1988). About the Prokhorov distance between the uniform distribution over the unit cube in $\rd$ and its empirical measure. \emph{Probab. Theory Related Fields} $\textbf{79}$, 431-450.

\bibitem{morris} \textsc{Morris, C.N.} (1983). Parametric empirical Bayes inference: theory and applications. \emph{J. Amer. Statist. Assoc.} $\mathbf{78}$, 47-65.

\bibitem{olkin} \textsc{Olkin, I.} and \textsc{Pukelsheim, F.} (1982). The distance between two random vectors with given dispersion matrices. \emph{Linear Algebra Appl.} $\mathbf{48}$, 257-263.

\bibitem{rizzelli} \textsc{Petrone, S., Rizzelli, S., Rousseau, J.} and \textsc{Scricciolo, C.} (2014). Empirical Bayes methods in classical and Bayesian inference. \emph{Metron} $\mathbf{72}$, 201-215.

\bibitem{petrone} \textsc{Petrone, S., Rousseau, J.} and \textsc{Scricciolo, C.} (2014). Bayes and empirical Bayes: do they merge? \emph{Biometrika} $\mathbf{101}$, 285-302.

\bibitem{robbins} \textsc{Robbins, H.} (1956). An empirical Bayes approach to statistics. \emph{Proc.} $3^{rd}$ \emph{Berkeley Symp.} $\mathbf{1}$, 157-163.

\bibitem{rockafellar} \textsc{Rockafellar, R.T.} (1970). \emph{Convex Analysis}. Princeton University Press, Princeton.

\bibitem{rousseau} \textsc{Rousseau, J.} (2016). On the frequentist properties of Bayesian nonparametric methods. \emph{Annual Review of Statistics and its Applications} $\mathbf{3}$, 211-231.

\bibitem{severini} \textsc{Severini, T.A.} (2000). \emph{Likelihood Methods in Statistics}. Oxford University Press, Oxford.

\bibitem{shorack} \textsc{Shorack, G.R.} and \textsc{Wellner, J.A.} (1986). \emph{Empirical Processes with Application to Statistics.} Wiley, New York.

\bibitem{shor} \textsc{Shor, P.W.} and \textsc{Yukich, J.E.} (1991). Minimax grid matching and empirical measures. \emph{Ann. Probab.} $\textbf{19}$, 1338-1348.

\bibitem{siegmund} \textsc{Siegmund, D.} (1969). On moments of the maximum of normed partial sums. \emph{Ann. Math. Statist.} $\textbf{40}$, 527-531.

\bibitem{talagrand} \textsc{Talagrand, M.} (1994). The transportation cost from the uniform measure to the empirical measure in dimension $\geq 3$. \emph{Ann. Probab.} $\textbf{22}$, 919-959.

\bibitem{talagrandIneq} \textsc{Talagrand, M.} (1996). Transportation cost for Gaussian and other product measures. \emph{Geom. Funct. Anal} $\textbf{6}$, 587-600.

\bibitem{teicher} \textsc{Teicher, H.} (1971). Completion of a dominated ergodic theorem. \emph{Ann. Math. Stat.} $\textbf{42}$, 2156-2158.

\bibitem{vdVaartWelln} van der \textsc{Vaart, A.W.} and \textsc{Wellner, J.A.} (1996). \emph{Weak Convergence and Empirical Processes.} Springer, New York.

\bibitem{VC} \textsc{Vapnik, V.} and \textsc{Chervonenkis, A.} (1971). On the uniform convergence of relative frequencies of events to their probabilities. \emph{Theory Probab. Appl.} $\textbf{16}$, 264-280.

\bibitem{villani} \textsc{Villani, C.} (2009). \emph{Optimal Transport. Old and New}. Springer, Berlin.

\bibitem{yosida} \textsc{Yosida, K.} and \textsc{Kakutani, S.}  (1939). Birkhoff's Ergodic Theorem and the Maximal Ergodic Theorem. \emph{Proc. Imp. Acad.} $\textbf{15}$, 165-168.

\bibitem{yukich} \textsc{Yukich, J.E.} (1989). Optimal matching and empirical measures. \emph{Proc. Amer. Math. Soc.} $\textbf{107}$, 1051-1059.

\bibitem{yurinskii} \textsc{Yurinskii, V.V.} (1976). Exponential inequalities for sums of random vectors. \emph{J. Multivariate Analysis} $\textbf{6}$, 473-499.

\bibitem{weedbach} \textsc{Weed, J.} and \textsc{Bach, F.}  (2017). Sharp asymptotic and finite-sample rates of convergence of empirical measures in Wasserstein distance. \emph{ArXiv:1707.00087}

\end{thebibliography}
\end{document}